\documentclass{amsart}
\usepackage{amsmath}
\usepackage{amssymb}
\usepackage{latexsym}
\usepackage{amsbsy}
\usepackage[all]{xy}
\usepackage{amscd}
\usepackage[dvips]{graphicx}
\newtheorem{theorem}{Theorem}[section]
\newtheorem{Lemma}[theorem]{Lemma}

\theoremstyle{definition}
\newtheorem{Definition}[theorem]{Definition}

\newtheorem{Prop}[theorem]{Proposition}
\newtheorem{Cor}[theorem]{Corollary}

\theoremstyle{remark}
\newtheorem{Remark}[theorem]{Remark}

\numberwithin{equation}{section}

\def\az{\alpha}      \def\ud{{\underline{d}}}
     \def\bd{{\bf d}}

\def\lz{\lambda}

\def\bbn{{\mathbb N}}  \def\bbz{{\mathbb Z}}  \def\bbq{{\mathbb Q}} \def\bb1{{\mathbb 1}}
\def\bbf{{\mathbb F}}   \def\bbe{{\mathbb E}} \def\bbf{{\mathbb F}}

  \def\bbc{{\mathbb C}}
  \def\fg{\mathfrak g} \def\fn{\mathfrak n}
\def\fh{\mathfrak h}

\def\ra{\rightarrow}
\def\hom{\mbox{Hom}}
\def\rad{\mbox{rad}\,} \def\ind{\mbox{ind}\,} 
   
\def\re{\mbox{\scriptsize re}}
 
\def\dim{\mbox{dim}\,}

\def\udim{\mbox{\underline {dim}}}

\def\ol{\overline}
\def\mod{\mbox{mod}\,}  
\def\Im{\mbox{Im}\,}    
\def\Ker{\mbox{Ker}\,}

\def\lr#1{\langle #1\rangle}

     \def\supp{\mbox{supp}\,}

\def\uq2{U_q(\hat{sl}_2)}

\def\bb{{\bf b}}

\def\bd{{\bf d}}
\def\nd{{\noindent}}

\def\blt{{\bullet}}
\def\pt{{\partial}}

\def\mc{{\mathcal{C}}}
\def\mk{{\mathcal{K}}}
\def\mn{{\mathcal{N}}}
\def\md{{\mathcal{D}}}
\def\mp{{\mathcal{P}}}
\def\mi{{\mathcal{I}}}
\def\mo{{\mathcal{O}}}
\def\mq{{\mathcal{Q}}}
\def\mr{{\mathcal{R}}}
\def\ml{{\mathcal{L}}}

\def\udd{{\underline{\bf d}}}
\def\ue{{\underline{e}}}
\def\uee{{\underline{\bf e}}}

\begin{document}

\title{Derived categories and Lie algebras}

\thanks{  The research was
supported in part by NSF of China and by the 973 Project of the
Ministry of Science and Technology of China.  }

\author{Jie Xiao, Fan Xu and Guanglian Zhang}
\address{Department of Mathematical Sciences\\
Tsinghua University\\
Beijing 10084, P.~R.~China} \email{jxiao@math.tsinghua.edu.cn
(J.Xiao),\  fanxu@mail.tsinghua.edu.cn (F.Xu)}
\address{Department of Mathematics\\
Shanghai Jiao Tong University\\
Shanghai 200240, P.~R.~China}
\email{zhangguanglian@mails.tsinghua.edu.cn (G.Zhang)}

\dedicatory{Dedicated to Professor George Lusztig }

\maketitle

\begin{abstract}
Let $\md^b(A)$ be the derived category of a finite dimensional basic
algebra $A$ with finite global dimension.  We construct the Lie
algebra arising from the 2-periodic version $\mk_2(\mp(A))$ of
$\mk^b(\mp(A))$ in term of constructible functions on varieties
attached to $\mk_2(\mp(A))$.
\end{abstract}


\bigskip

\section{ Introduction}
\subsection{}
In the last thirty years of the twentieth century, there were two
parallel fields in mathematics got extensively developed. One is the
infinite dimensional Lie theory, in particular, the Kac-Moody Lie
algebras. One is the representation theory of finite dimensional
algebras, in particular, the representations of quivers. The close
relation between the two subjects was discovered in a very early
stage. Gabriel in \cite{Gabriel1972} found that the quivers of
finite representation type were given by the Dynkin graphs in Lie
theory, and the dimension vectors provide the bijective
correspondence between the isomorphism classes of indecomposable
representations of the quiver and the positive root system of the
semisimple Lie algebra. After the Gabriel theorem, a lot of progress
on the connection between the representations of quivers or
hereditary algebras and Lie algebras had been made, for example, by
Bernstein-Gelfand-Ponomarev \cite{BGP1973}  and Dlab-Ringel
\cite{DlabRingel}. The final and most general result is the Kac
theorem \cite{kac1980} which extends to consider the quiver and the
symmetric Kac-Moody algebra of arbitrary type. It states that the
dimension vectors of indecomposable representations are exactly the
positive roots, a unique indecomposable corresponds to each real
root and infinitely many to each imaginary root; the multiplicity of
the imaginary root, which is conjectured by Kac in \cite{kac1982},
is given by a geometric parameter in terms of representations of the
quiver. A new progress on the Kac conjecture is by Crawley-Boevey
and Van den Bergh in \cite{CrawleyBoevey}.

\subsection{}
Ringel in \cite{Ringel1990} discovered his Hall algebra structure
by giving an answer to the following fundamental question: how to
recover the underlying Lie algebra structure directly from the
category of representations of the quiver.

Let $Q$ be a quiver, $A=\mathbb{F}_qQ$ the path algebra of $Q$
over $\mathbb{F}_q:$ the finite field with $q$ elements. Set
$\mathcal{P}=\{{\rm isoclasses \ of\ representations\ of}\ Q\}.$
For any $\alpha\in\mathcal{P}$ choose $V_\alpha$ to be a
representative in the class $\alpha.$ Given three classes
$\lambda,\ \alpha,\ \beta\in\mathcal{P},$ let
$g^\lambda_{\alpha\beta}$ be the order of the finite set $\{ W
\lhd V_\lambda|W\cong V_\beta,V_\lambda/W\cong V_\alpha\}.$ By
taking $v=\sqrt{q}$ and the integral domain $\mathbb{Q}(v),$ the
(twisted) Ringel-Hall algebra $\mathcal{H}^*(A)$ can be defined to
be a free $\mathbb{Q}(v)$-module with basis
$\{u_\lambda|\lambda\in \mathcal{P}\}$ and multiplication is given
by
$$u_\alpha *
u_\beta=v^{\lr{\alpha,\beta}}\sum_{\lambda\in
\mathcal{P}}g^\lambda_{\alpha\beta}u_\lambda\text{ for \ all }
\alpha,\beta\in \mathcal{P}.$$

  One may consider the subalgebra of $\mathcal{H}^*(A)$ generated
by $u_i=u_{\alpha_i},$ for $i\in I (=Q_0)��$ where
$\alpha_i\in\mathcal{P}$ is the isoclass of simple $A$-module at
vertex $i.$ The subalgebra is called the {\em composition algebra}
and it is denoted by  $\mathcal{C}^\ast(A).$ On the other hand,
the index set $I$ of simple $A$-modules together with the
symmetric Euler form $(I,(-,-))$ of $A$ is a Cartan datum in the
sense of Lusztig \cite{Lusztig1993}. For a Cartan datum
$(I,(-,-)),$ the quantized enveloping algebra $\mathcal{U}_q$
defined by Drinfeld \cite{Drinfeld1988} and Jimbo \cite{Jimbo1985}
is associated with it. The positive part $\mathcal{U}_q^+$ is
generated by $E_i,$ $i\in I$ with subject to the quantum Serre
relations.

There is a usual way to define the generic form
$\mathcal{C}^\ast(Q)$ of the composition algebra
$\mathcal{C}^\ast(A)$ by considering the representations of $Q$
over infinitely many finite fields. Then $\mathcal{C}^\ast(Q)$ is
a $\mathbb{Q}(v)$-algebra where $v$ becomes a transcendental
element over $\mathbb{Q}.$  Put
$u_i^{(*n)}=\frac{u_i^{*n}}{[n]_i!}$ for $i\in I$ and
$n\in\mathbb{N}$ and let $\mathcal{C}^\ast(Q)_{\mathcal{Z}}$ be
the integral form of $\mathcal{C}^\ast(Q),$ which is generated by
$u_i^{(*n)},\ i\in I,\ n\in\mathbb{N}$ over the integral domain
$\mathcal{Z}=\mathbb{Z}[v,v^{-1}].$ Also the quantum group
$\mathcal{U}_q^+$ has the integral form
$\mathcal{U}^+_{\mathcal{Z}},$ which is generated by $E_i^{(n)},\
i\in I,\ n\in\mathbb{N}$ over $\mathcal{Z}.$ Then by Ringel
\cite{Ringel1990a} and Green \cite{Green1995}, the canonical map
$\mathcal{C}^{*}(Q)_{\mathcal{Z}}\rightarrow
\mathcal{U}^+_{\mathcal{Z}}$ by sending $u_i^{(*n)}$ to $
E_i^{(n)}$ for $i\in I$ and $n\in\mathbb{N}$ leads to  a
$\mathcal{Z}$-algebra isomorphism, if the two algebras share a
common Cartan datum.

Let $\ind\mp=\{{\rm isoclasses \ of\ indecomposable\
representations\ of}\ Q\}.$ Then $g^\lambda_{\alpha\beta}$ can be
regarded as a function on $q$ for $\alpha, \beta,
\lambda\in\ind\mp.$ In fact, Ringel in \cite{Ringel1990} proved
that $g^\lambda_{\alpha\beta}$ is an integral polynomial on $q$
when $Q$ is of finite type. One can take the integral value
$g^\lambda_{\alpha\beta}(1)$ by letting that $q$ tends to $1.$
Ringel \cite{Ringel1990}, for $Q$ of finite type, prove that the
$u_{\alpha}, \alpha\in\ind\mp,$ spanned a Lie subalgebra of
$\mathcal{C}^\ast(Q)\mid_{q=1}$ with Lie bracket

$$[u_{\alpha},u_{\beta}]=\sum_{\lambda\in\ind\mp}(g^\lambda_{\alpha\beta}(1)-g^\lambda_{\beta\alpha}(1))u_{\lambda}$$
for $\alpha, \beta\in\ind\mp.$ This realized the positive part
$\mathfrak{n^+}$ of the semisimple Lie algebra $\mathfrak{g}.$ Of
course, Ringel's approach also works for $Q$ of arbitrary type. In
general, there exists the generic composition Lie subalgebra
$\mathcal{L}$ of $\mathcal{C}^\ast(Q)\mid_{q=1}$ generated by
$u_i, i\in I$ and $\ind\mp$ is no longer to index a basis of
$\mathcal{L}.$ Now $\mathcal{L}$ is canonically isomorphic to the
positive part $\mathfrak{n^+}$ of the symmetric Kac-Moody Lie
algebra $\mathfrak{g}.$ For a realization of the whole
$\mathfrak{g},$ not just its positive part, Peng and Xiao in
\cite{PengXiao2000} have constructed a Lie algebra from a
triangulated category with the 2-periodic shift functor $T,$ i.e,
$T^2=1.$ If, specially, consider the 2-periodic orbit category of
the derived category of a finite dimensional hereditary algebra,
the Lie algebra obtained in \cite{PengXiao2000} gives rise to the
global realization of symmetrizable Kac-Moody algebra of arbitrary
type. In \cite{PengXiao2000} they consider the triangulated
categories over finite fields. Replacing counting the order of the
filtration set, they calculate the order of the orbit space of a
triangle. By a hard work, they obtained a Lie ring
$\mathfrak{g}_{(q-1)}$ over $\bbz/(q-1)$ for the prime powers
$q=|\bbf_q|.$ Then they performed their work over finite field
extensions of arbitrarily large order and construct a generic Lie
algebra which is similar to the generic composition Lie subalgebra
done by Ringel in \cite{Ringel1992}. A transcendental Lie algebra
was finally obtained.

\subsection{}
Quickly after the work of Ringel \cite{Ringel1990a}, people
realized that a geometric setting of Ringel-Hall algebra is
possible by using the convolution multiplication (see
\cite{Schofield} and \cite{lusztig1990}). Let $Q$ be a quiver and
$\alpha=\sum_{i\in I}a_i i\in\mathbb{N}[I]$  a dimension vector.
We fix a $I$-graded space
$\mathbb{C}^{\alpha}=(\mathbb{C}^{a_i})_{i\in I}.$ Then
$$\bbe_\alpha=\bigoplus_{h:s(h)\rightarrow t(h)}{\rm Hom}_{\mathbb{C}}
(\mathbb{C}^{a_{s(h)}},\mathbb{C}^{a_{t(h)}})$$ is an affine
space. Set
$$G_\alpha=\Pi_{i\in I}GL(a_i,\mathbb{C}).$$
For any $(x_h)\in \bbe_\alpha$ and $g=(g_i)\in G_{\alpha},$ we
define the action $g\cdot(x_h)=(g_{t(h)}x_hg^{-1}_{s(h)}).$ For
any $Q$-representation $M$ with $\udim M=\alpha,$ let
$\mathcal{O}_M\subset\bbe_\alpha$ be the $G_\alpha$-orbit of $M.$

For an algebraic variety $X$ over $\mathbb{C},$ a subset $A$ of
$X$ is said to be constructible if it is a finite union of locally
closed subsets. A function $f: X\rightarrow \mathbb{C}$ is
constructible if it is a finite $\mathbb{C}$-linear combination of
characteristic functions $\mathbf{1}_{\mathcal{O}}$ for
constructible subsets $\mathcal{O}.$

We define $M_{G_\alpha}(Q)$ to be the space of constructible
$G_\alpha$-invariant functions $\bbe_\alpha\rightarrow \mathbb{C}
,$ and let $M_G(Q)=\bigoplus_{\alpha\in
\mathbb{N}I}M_{G_\alpha}(Q).$ Let ${\rm ind}\bbe_\alpha(Q)$ to be
the constructible subset of $\bbe_\alpha$ consisting of all points
$x$ which correspond to indecomposable $Q$-representations, and
let ${\rm ind}M_{G_\alpha}(Q)$ to be the space of constructible
$G_\alpha$-invariant functions over ${\rm ind}\bbe_\alpha.$ We may
regard as ${\rm ind}M_{G_\alpha}(Q)=\{f\in M_{G_\alpha}(Q)|{\rm
supp}f\subseteq {\rm ind}\bbe_\alpha\},$ and ${\rm ind}
M_G(Q)=\oplus_{\alpha\in R^+}{\rm ind} M_{G_\alpha}(Q), $ where,
by Kac theorem, $R^+$ is the positive root system of the Kac-Moody
Lie algebra corresponding to $Q.$ The space
$M_G(Q)=\bigoplus_{\alpha\in \mathbb{N}I}M_{G_\alpha}(Q)$ cab be
endowed with the associative algebra structure by the convolution
multiplication:

$$1_{\mathcal{O}_1}*1_{\mathcal{O}_2}(y)=\chi(\mathcal{F}_{\mo_1
\mo_2}^{y})$$ for any $G_\alpha$-invariant constructible set
$\mo_1$ and $G_\beta$-invariant constructible set $\mo_2$ with
$\alpha,\beta\in\bbn I,$ where $\mathcal{F}_{\mo_1
\mo_2}^{y}=\{x\in\mo_2 \mid M(x)\subseteq M(y)\ \mbox{and}\
M(y)/M(x)\in \mo_1\}$ and $\chi(X)$ denotes the Euler
characteristic of the topological space $X.$ As in
\cite{Riedtmann1994} and \cite{dxx2006}, it can be proved that the
space ${\rm ind} M_G(Q)$ has a Lie algebra structure under the
usual Lie bracket
$$[1_{\mo_1},1_{\mo_2}]=1_{\mo_1}*1_{\mo_2}-1_{\mo_2}*1_{\mo_1}.$$
Applying this setup to the case $Q$ being a tame quiver,
Frenkel-Malkin-Vybornov \cite{FMV2001} gave an explicit
realization of the positive parts of affine Lie algebras.

\subsection{}
The great progress is made by Lusztig, who apply the Hall algebras
in a geometric setting to study the quantum groups (see
\cite{lusztig1990} and \cite{Lusztig1991a}) and the enveloping
algebras (see \cite{lusztig2000}). The canonical bases of the
quantum groups and the semicanonical bases of the enveloping
algebras were originally constructed in terms of representations
of quivers. However Lusztig \cite{Lusztig1991a} has pointed out
that a more suitable choice is the preprojective algebra, which is
given by the double quiver of $Q$ with the Gelfand-Ponomarev
relations. Further progress in this direction is the study of
Nakajima \cite{Nakajima1998} on his quiver varieties, which leads
to a geometric realization of the representation theory of
Kac-Moody algebras.

Inspired by  Ringel's work on Hall algebras and Lusztig's geometric
approach to quantum groups, the aim of this paper is to give a
global and  geometric realization of the Lie algebras arising from
the derived categories, which is a generalization of our earlier
work \cite{PengXiao2000}.

\subsection{}

If we consider the module category of $A=\bbc Q/J,$ we have the
algebraic variety $\bbe_{\ud}(Q,R)$ for $A$-modules with a fixed
dimension vector $\ud$ and it is a $G$-variety where $G=G_{\ud}(Q)$
is a reductive group. According to the work of C.de Concini and
E.Strickland in \cite{CS1981} and M.Saorin and B.Huisgen-Zimmermann
in \cite{SHZ2001}, this geometry can be generalized to over the
chain complexes of $A$-modules. Section 2 is devoted to do this. Let
$K_0(\md^b(A))$ be the Grothendieck group of $\md^b(A)$ and $\udim$
be the canonical map from the abelian group of dimension vector
sequences to $K_0(\md^b(A))$. Given $\bd\in K_0(\md^b(A))$ and
$\udd\in\udim^{-1}(\bd)$, the set $\mc^b(A,\udd)$ of all complexes
of $A$-modules with the dimension vector sequences $\udd$ and its
subset $\mp^b(A,\udd)$ of all projective complexes  can be endowed
with the affine variety structures. Let $\mc^b(A, \bd)$ (resp.
$\mp^b(A, \bd)$) be the direct limit of $\mc^b(A, \udd)$ (resp.
$\mp^b(A, \udd)$) for $\udd\in \udim^{-1}(\bd).$ Here, we associate
to $\mk^b(\mp(A))$ its quotient space $\mq\mp^b(A,\bd)$ which is the
direct limit of the quotient spaces $\mq\mp^b(A, \udd))$ of
$\mp^b(A, \udd)$ under the action of some algebraic group
$G_{\udd}$. Our main aims in Section 2 are  to study the relation
between $\mq^b(A,\bd)$ and $\mq\mp^b(A,\bd)$,  the action of derived
equivalence on $\mq^b(A,\bd)$ and $\mq\mp^b(A,\bd)$ and characterize
the orbit space $\mq\mp^b(A,\bd)$ of $\mp^b(A, \bd)$ under the
action of the direct limit $G_{\bd}$ of algebraic groups $G_{\udd}$.
Therefore the main point in Section 2 is that, we can regard the
$G_{\bd}$-invariant geometry in $\mp^b(A,\bd)$ as the moduli space
in which the orbits index the isomorphism classes of objects in the
derived category.

In Section 3, we consider the inverse limit of the $\bbc$-space of
$G_{\uee}$-invariant constructible functions over $\mp^b(A,\uee)$
for any $\uee\in \udim^{-1}(\bd).$  Any element in the inverse limit
can be viewed as a $G_{\bd}$-invariant constructible function over
$\mp^b(A,\bd)$. We define the convolution between
$G_{\bd}$-invariant constructible functions.  Our main result in
Section 3 is that our convolution rule is well-defined. In order to
prove it, we need to define the \emph{na\"ive} Euler characteristic
of the orbit spaces induced by triangles in the triangulated
category as in \cite[Section 4.3]{Joyce}. The theorem of Rosenlicht
\cite{Rosenlicht1963} for the algebraic group action on varieties is
crucial for us.  Section 4 is just to transfer the results in
Section 3 to the 2-periodic orbit categories of the derived
categories.

Section 5 is devoted to verifying the Jacobi identity. In
\cite{PengXiao2000}, by counting the Hall numbers $F_{XY}^L$ for the
triangles of the form $X\rightarrow L\rightarrow Y\rightarrow X[1]$,
it has been proved that the Jacobi identity can be deduced from the
octahedral axiom of the triangulated categories. However we need to
prove that the correspondences among the various orbit spaces in the
derived categories are actually given by the algebraic morphisms of
algebraic varieties. We think this geometric method is more
transparent to reflect the hidden symmetry in the derived category.
Additionally, we get the two properties which is unknown
 in \cite{PengXiao2000}. Firstly the proper assumption in \cite{PengXiao2000} is not necessary, in
 fact, it is easy to give examples such that $\udim X=0$ for some
 nonzero indecomposable $X$ in $\md^b(A).$ Secondly we show that
 the Lie algebras arising from the 2-periodic orbit categories of
 the derived categories always possess the symmetric invariant
 form in the sense of Kac \cite{kac1990}, which is essentially
 non-degenerated. Section 6 is to apply the construction to
 the 2-periodic orbit categories of the derived categories of
 representations of quivers, particularly, tame quivers. This
 gives rise to a global realization of the symmetric (generalized) Kac-Moody algebras of arbitrary type.
 In particular, an explicit realization of the affine Lie
 algebras.
\subsection{}

 Finally we should mention recent advances by To\"en \cite{Toen2005} and
 Joyce \cite{joyce2003}. To\"en defined  an associative algebra, called the derived Hall algebra
 associated to a dg category over a finite field. A direct proof for To\"en's theorem is given in \cite{XX2006}.
 Joyce considered a new Ringel-Hall type
 algebra consisting of functions over stacks associated to abelian
 categories. Their results can be viewed as  improvements of the Ringel-Hall type algebra
 with respect to categorification and geometrization.
However, it is unknown how to define an analogue of the derived Hall
over the complex field (see \cite{Lusztig2000-02}
\cite{Nakajima1998}) or an analogue of the derived algebra for the
2-period version of a derived category (see \cite{Toen2005} and
\cite{XX2006}). Hence, it is still an open question to define an
associative multiplication which induces the Lie bracket in this
paper and supplies the realization of the corresponding enveloping
algebra.

\section{Topological spaces attached to derived categories}

\subsection{Module varieties}
Given an associative algebra $A$ over the complex field $\bbc$ ,
in this paper, we always assume that $A$ is both finite
dimensional and finite global dimensional. By a result of
P.Gabriel (\cite{Gabriel1972}) the algebra $A$ is given by a
 quiver $Q$ with relations $R$ (up to Morita equivalence). Let $Q=(Q_0,Q_1,s,t)$ be a quiver, where
$Q_0$ and $Q_1$ are the sets of vertices and arrows, respectively, and $s,t: Q_1\rightarrow Q_0$ are
 maps such that any arrow $\az$ starts at $s(\az)$ and terminates at $t(\az).$ For any dimension vector
 $\ud=(d_i)_{i\in Q_0},$ we consider the affine space over $\bbc$
$$\bbe_{\ud}(Q)=\bigoplus_{\az\in Q_1}\hom_{\bbc}(\bbc^{d_{s(\az)}},\bbc^{d_{t(\az)}})$$
Any element $x=(x_{\az})_{\az\in Q_1}$ in $\bbe_{\ud}(Q)$ defines
a representation $(\bbc^{\ud}, x)$ where
$\bbc^{\ud}=\bigoplus_{i\in Q_0}\bbc^{d_i}$. A relation in $Q$ is
a linear combination $\sum_{i=1}^{r}\lz_{i}p_i,$ where
$\lz_i\in\bbc$ and $p_i$ are paths of length at least two with
$s(p_i)=s(p_j)$ and $t(p_i)=t(p_j)$ for all $1\leq i,j\leq r.$ For
any $x=(x_{\az})_{\az\in Q_1}\in\bbe_{\ud}$ and any path
$p=\az_1\az_2\cdots\az_m$ in $Q$ we set
$x_{p}=x_{\az_1}x_{\az_2}\cdots x_{\az_m}.$ Then $x$ satisfies a
relation
 $\sum_{i=1}^{r}\lz_{i}p_i$ if $\sum_{i=1}^{r}\lz_i x_{p_i}=0.$ If $R$ is a set of relations in $Q,$ then
let $\bbe_{\ud}(Q,R)$ be the closed subvariety of $\bbe_{\ud}(Q)$ which consists of all elements satisfying
all relations in $R.$ Any element $x=(x_{\az})_{\az\in Q_1}$ in $\bbe_{\ud}(Q,R)$ defines in a natural way a
representation $M(x)$ of $A=\bbc Q/J$ with $\udim M(x)=\ud,$ where $J$ is the admissible ideal generated by $R.$
 We consider the algebraic group
$$G_{\ud}(Q)=\prod_{\i\in Q_0}GL(d_i,\bbc),$$ which acts on
 $\bbe_{\ud}(Q)$ by $(x_{\az})^{g}=(g_{t(\az)}x_{\az}g_{s(\az)}^{-1})$ for $g\in G_{\ud}$ and $(x_{\az})\in\bbe_{\ud}.$
It naturally induces the action of $G_{\ud}(Q)$ on
$\bbe_{\ud}(Q,R).$ The induced orbit space is denoted by
$\bbe_{\ud}(Q,R)/G_{\ud}(Q).$ There is a natural bijection between
the set ${\mathcal M}(A,\ud)$ of isomorphism classes of
$\bbc$-representations of $A$ with dimension vector $\ud$ and the
set of orbits of $G_{\ud}(Q)$ in $\bbe_{\ud}(Q,R).$ So we may
identify  ${\mathcal M}(A,\ud)$ with $\bbe_{\ud}(Q,R)/G_{\ud}(Q).$

\subsection{Categories of complexes}
First we consider the category of complexes $\mathcal{C}(A)$.
 Its objects are sequences $M^{\bullet}=(M_n,\pt_n)$ of finite dimensional $A$-modules and their homomorphisms
 \begin{equation}\label{e11}
  \begin{CD}
  \dots  @>\pt_{n-1}>> M_{n} @>\pt_{n}>> M_{n+1} @>\pt_{n+1}>> M_{n+2} @>\pt_{n+2}>> \dots
 \end{CD}
 \end{equation}
 such that $\pt_{n+1}\pt_n=0$ for all $n$. A morphism $\phi^{\blt}:M^{\blt}\ra M'^{\blt}$ between two complexes is a sequence of homomorphisms $\phi^{\blt} =(\phi_n:M_n\ra M'_n)_{n\in\bbz}$ such that the following diagram is commutative.

 \begin{equation}\label{e12}
  \begin{CD}
  \dots  @>\pt_{n-1}>> M_{n} @>\pt_{n}>> M_{n+1} @>\pt_{n+1}>> M_{n+2} @>\pt_{n+2}>> \dots  \\
  &&    @V\phi_{n}VV  @V\phi_{n+1}VV  @V\phi_{n+2}VV \\
  \dots  @>\pt'_{n-1}>> M'_{n} @>\pt'_{n}>> M'_{n+1} @>\pt'_{n+1}>> M'_{n+2} @>\pt'_{n+2}>> \dots
 \end{CD}
 \end{equation}
 One says that such a morphism is \emph{homotopic to zero} if there are homomorphisms
 $\sigma_n:M_n\to M'_{n-1} $  such that $\phi_n=\sigma_{n+1}\pt_n+\pt'_{n-1}\sigma_n$ for all
 $n\in\bbz$. The factor category $\mathcal{K}(A)$ of $\mathcal{C}(A)$ modulo the ideal of morphisms homotopic to zero
 is called the \emph{homotopic category} of $A$-modules. For each $n$ the \emph{$n$-th homology}
 of a complex is defined as $H_n(M^{\blt})=\Ker \pt_n/\Im \pt_{n-1}$. Obviously,
 a morphism $\phi^{\blt}$ of complexes induces homomorphisms of homologies $H_n(\phi^{\blt}):
 H_n(M^{\blt})\rightarrow H_n(M'^{\blt})$ and if $\phi^{\blt}$ is homotopic to zero, it induces zero
homomorphisms of homologies. One call a morphism $\phi^{\blt}$
in $\mc(A)$
 or in $\mk(A)$ \emph{quasi-isomorphism} if the induced morphisms $H_n(\phi^{\blt})$ are isomorphisms for all $n.$
 Now the derived category $\md(A)$ is defined to be the category of fractions $\mk(A)[\mn^{-1}]$,
 where $\mn$ is the set of all quasi-isomorphisms, which is obtained from $\mk(A)$ by inversing all morphisms in $\mn.$
 One calls a complex \emph{right bounded} (\emph{left bounded}, \emph{bounded}, respectively) if there is $n_0$ such
that $M_n=0$ for $n>n_0$ (
there is $n_1$ such that $M_n=0$ for
 $n<n_1$, or there are both, respectively). The corresponding categories are denoted by $\mc^-(A),\mk^-(A),
 \md^-(A)$ ( by $\mc^+(A),\mk^+(A),\md^+(A)$, or by $\mc^b(A),\mk^b(A),\md^b(A),$ respectively). In this paper we
 mainly deal with the bounded situation.

 The category $A$-mod of finite dimensional $A$-modules can be naturally embedded into $\md(A)$ (even
 in $\md^b(A)$): a module $M$ is identified with the complex $M^{\blt}$ such that $M_0=M$ and $
 M_n=0$ for $n\ne0$.

 A complex $P^{\blt}=(P_n,\pt_n)$ is called projective if all $P_n$ are projective $A$-modules.
 Since the category $A$-mod has enough projective objects, one can replace, when considering right bounded
 homotopic and derived category, arbitrary complexes by projective ones.
 We denote by $\mp^-(A)$ and by $\mp^b(A)$ the full subcategories of $\mc^-(A)$ and $\mc^b(A)$ which consist of
 right bounded and bounded
 projective complexes, respectively. Actually, we have $\md^-(A)\simeq \mk^-(\mp^-(A))\simeq \mp^-(A)/\mi$,
 where $\mi$ is the
 ideal of morphisms homotopic to zero (see\cite{GM1996}). Moreover, every  finite dimensional $A$-module
 $M$ has a \emph{projective cover}, i.e., an epimorphism $p_M:P(M)\to M$  such that  $P(M)$
 is projective and $\Ker p_M\subseteq\rad P(M)$, the radical of $P(M)$. Therefore, we
 can only consider \emph{minimal  or radical }  projective complexes $P^{\blt}=(P_n,\pt_n)$ with the property:
 $P_n$ is projective and  $\Im \pt_n\subseteq\rad
 P_{n+1}$ for all $n$. Let  $\rad\mp^-(A)$ be the full subcategory of $\mp^-(A)$ which consist of minimal
 projective complexes. Since every projective complex in $\mp^-(A)$ is quasi-isomorphic to a minimal projective
 complexes, we have $\md^-(A)\simeq \rad\mp^-(A)/\mi,$ where $\mi$ is the
 ideal of morphisms homotopic to zero.  One immediately checks that  a
 morphism $\phi^{\blt}$ between  minimal projective complexes induces an isomorphism in
 $\md^-(A)$ if and only if
 $\phi^{\blt}$ itself is an isomorphism in $\rad\mp^-(A).$ If we further assume that the global dimension of $A$ is
 finite, then we have
 $\md^b(A)\simeq \rad\mp^b(A)/\mi$, since any bounded
 complex has a bounded projective resolution.
\\

\subsection{Complex varieties}    The geometrization of $A$-modules with dimension vector $\ud$ can be carried over, in the same spirit,
to the complexes. Let the algebra $A=\bbc Q/J$ be of finite global
dimension and the admissible ideal $J$  is given by a set $R$ of
relations in $Q.$ For a dimension vector $\ud$ we understand as
$\ud:Q_0\ra\bbn.$ We set the $Q_0$-graded $\bbc$-space
$\bbc^{\ud}=\bigoplus_{j\in Q_0}\bbc^{\ud(j)}.$ For a sequence of
dimension vectors $\udd=(\cdots,\ud_{-1},\ud_{0},\ud_1,\cdots)$ with
only finite many non-zero entries, we define $\mc^b(A,\udd)$ to be
the closed subset of (see \cite{SHZ2001})
$$\prod_{i\in\bbz}\bbe_{\ud_i}(Q,R)\times\prod_{i\in\bbz}\hom_{\bbc}(V^{\ud_i},V^{\ud_{i+1}})$$
which consists of elements $(x_i,\pt_i)_i,$ where
$x_i\in\bbe_{\ud_i}(Q,R)$ and $M(x_i)=(\bbc^{\ud_i}, x_i)$ is the
corresponding $A$-module and
$\pt_i\in\hom_{\bbc}(\bbc^{\ud_i},\bbc^{\ud_{i+1}})$ is a
$A$-module homomorphism from $M(x_i)$ to $M(x_{i+1})$ with the
property $\pt_{i+1}\pt_{i}=0.$ In fact, $(M(x_i),\pt_i)_i,$ or
simply denoted by $(x_i,\pt_i)_i,$ is a complex of $A$-modules and
$\udd$ is called its dimension vector sequence.

The group $G_{\udd}:=\prod_{i\in\bbz} G_{\ud_i}(Q)$ acts on
$\mc^b(A,\udd) $ via the conjugation action
$$(g_i)_i (x_i,\pt_i)_i=((x_i)^{g_i},g_{i+1}\pt_i g_i^{-1})_i$$
where the action $(x_i)^{g_i}$ was defined as in Section 2.1.
Therefore the orbits under the action correspond bijectively to
the isomorphism classes of complexes of $A$-modules.

We fix a set $P_1,P_2,\cdots,P_l$ to be a complete set of
indecomposable projective $A$-modules (up to isomorphism). Let
$\mp^b(A)$ be the full subcategory of $\mc^b(A)$ which consists of
projective complexes $P^{\blt}=(P^i,\pt_i)$ such that each $P^i$ has
the decomposition $P^i\cong\bigoplus_{j=1}^l e_j^i P_j.$ We denote
by $\ue(P^i)$ the vector $(e^i_1,e^i_2,\cdots,e^i_l).$ The sequence
$(\cdots,\ue(P^{-1}),\ue(P^0),\ue(P^1),\cdots)$, denoted by
$\uee(P^{\blt}),$
  is called the projective dimension sequence of
$P^{\blt}.$ Put $\udd(\uee)=(\ud_i),$ where $\ud_i=\udim P^i,$. In
the similar way as in \cite{JSZ2005}, for a fixed projective
dimension sequence $\uee=(\cdots, \ue_i, \cdots),$ we define
$\mp^b(A,\uee)$ to be the locally closed subset of
$\mc^b(A,\udd(\uee))$ consisting of $(x_i,\pt_i)_i$ with
$(\bbc^{\ud_i}, x_i)$ isomorphic to $P^i$ for any $i\in \bbz.$  The
action of the algebraic group $G_{\udd(\uee)}:=\prod_{i\in
\bbz}G_{\ud_i}$ on $\mc^b(A,\udd(\uee))$ induces an action on
$\mp^b(A,\uee).$

\subsection{} In this subsection we consider the topological structures which are endowed with
  $\mc^b(A)$ and $\mp^b(A)$.

Let $K_0(\md^b(A)),$ or simply by $K_0,$ be the Grothendieck group
of the derived category $\md^b(A),$ and $\udim:\md^b(A)\ra
K_0(\md^b(A))$ the canonical surjection. It induces a canonical
surjection from the abelian group of dimension vector sequences to
$K_0,$ we still denote it by $\udim.$ We have known that the set
$\mc^b(A,\udd)$ of all complex of fixed dimension vector sequence
$\udd$ in $\mc^b(A)$ is an affine variety.  For any
$\udd_1,\udd_2\in \udim^{-1}(\bd),$ we write $\udd_1\leq\udd_2,$ if
there exists a complex $M^{\bullet}(\udd_1,\udd_2)$ in $\mc^b(A,
\udd_2-\udd_1)$ which is a direct sum of shifted copies of complexes
of the form
$$\xymatrix{\cdots\ar[r]&0\ar[r]&S\ar[r]^{1}&S\ar[r]&0\ar[r]&\cdots}
$$
where $S$ is a simple $A$-module. This defines a partial order on
$\udim^{-1}(\bd)$. Fix the set $$\{M^{\bullet}(\udd_1,\udd_2)\mid
\udd_1,\udd_2\in \udim^{-1}(\bd), M^{\bullet}(\udd_1,\udd_2)\oplus
M^{\bullet}(\udd_2,\udd_3)=M^{\bullet}(\udd_1,\udd_3)\}.$$ We have a
morphism of varieties :
$$ T_{\udd_1\udd_2}: \mc^b(A, \udd_1)\rightarrow \mc^b(A,\udd_2)$$
mapping a complex $X^{\bullet}$ to $X^{\bullet}\oplus M^{\bullet}(\udd_1,\udd_2)$.
Then we obtain a direct system $\{(\mc^b(A,\udd), T_{\udd\udd'})\mid
\udd,\udd'\in \ \udim^{-1}(\bd)\}$ and define
$$\mc^b(A,\bd)=\varinjlim_{\udd\in\udim^{-1}(\bd)}\mc^b(A,\udd)$$ for $\bd\in K_0$.
We have a canonical morphism $T_{\udd}: \mc^b(A, \udd)\rightarrow
\mc^b(A, \bd)$ for any $\udd\in \udim^{-1}(\bd).$ A subset $U$ is
open in $\mc^{b}(A, \bd)$ if and only if $T_{\udd}^{-1}(U)$ is open
in $\mc^b(A, \udd)$ for any $\udd\in \udim^{-1}(\bd).$

  Moreover, we also define the quotient space $\mq^b(A,\bd)=\mc^b(A,\bd)/\sim,$ where $x\sim y$ if and only if the
corresponding complexes $M(x)^{\blt}$ and $M(y)^{\blt}$ are
quasi-isomorphic to each other, i.e., they are isomorphic in
$\md^b(A).$ The topology of $\mq^b(A,\bd)$ is quotient topology,
i.e., let $\pi: \mc^b(A,\bd)\rightarrow \mq^b(A,\bd)$ be the
canonical surjection, $U$ is an open (closed) set of $\mq^b(A,\bd)$
if and only if $\pi^{-1}(U)$ is an open (closed) set of
$\mc^b(A,\bd)$. Let $S$ be a subset of $\mc^b(A,\bd)$. Its orbit
space is defined as $\mo(S)=\{y|y\sim x\mbox{ for some}\ x\in S \}$.

A complex $X^{\blt}=(X_i,\pt_i)_i$ is called contractible if the
induced homological groups $H_i(X^{\blt})=\ker\pt_{i+1}/ \Im\pt_i=0$
for all $i\in\bbz.$ It is easy to see that any contractible
projective complex is isomorphic to a direct sum of shifted copies
of complexes of the form
\begin{equation}\nonumber\label{e11}
  \begin{CD}
  \dots  @> >>0 @>0>> P @>f>> P @>0 >>0 @> >>\dots
 \end{CD}
 \end{equation}
with $P$ is a projective $A$-module and $f$ is an automorphism. We
call an element $x\in\mp^b(A,\uee)$ contractible if the
corresponding projective complex is contractible. Let $\uee_1$ and
$\uee_2$ be projective dimension sequence such that $\udd(\uee_1)$
and $\udd(\uee_2)$ are in $\udim^{-1}(\bd)$. We write
$\uee_1\leq\uee_2$ if there exists a contractible projective complex
$P^{\bullet}(\uee_1,\uee_2)\in \mp^b(A, \uee_2-\uee_1).$ We fix the
set $$\{P^{\bullet}(\uee_1, \uee_2)\mid P^{\bullet}(\uee_1,
\uee_2)\oplus P^{\bullet}(\uee_2, \uee_3)=P^{\bullet}(\uee_1,
\uee_3), \udd(\uee_1), \udd(\uee_2),\udd(\uee_3)\in
\udim^{-1}(\bd)\}$$ of contractible projective complexes. Then we
have a canonical morphism of varieties
$$t_{\uee_1\uee_2}: \mp^b(A, \uee_1)\rightarrow \mp^b(A, \uee_2)$$
mapping $X^{\bullet}$ to $X^{\bullet}\oplus P^{\bullet}(\uee_1,
\uee_2).$  Hence, we can define
$$\mp^b(A,\bd)=\varinjlim_{\udd(\uee)\in\udim^{-1}(\bd)}\mp^b(A,\uee)$$
for any $\bd\in K_0.$  By definition, there are canonical morphisms
$t_{\uee}: \mp^b(A, \uee)\rightarrow \mp^b(A, \bd).$ We have the
quotient space
$$\mq\mp^b(A,\bd)=\mp^b(A,\bd)/\sim,$$
where $x\sim y$ in $\mp^b(A,\bd)$ if and only if the corresponding projective complexes $P(x)^{\blt}$ and $P(y)^{\blt}$
are quasi-isomorphic, i.e., they are isomorphic in $\md^b(A).$ The topology for $\mq\mp^b(A,\bd)$ is the quotient
topology from $\mp^b(A,\bd).$ For any $x\in\mp^b(A,\bd),$ then the orbit is
$$\mo(x)=\{y\in\mp^b(A,\bd)|y\sim x\}.$$

For $\uee_1\leq\uee_2$, there exist a natural morphism of varieties
$f_{\uee_1\uee_2}: G_{\udd(\uee_1)}\rightarrow G_{\udd(\uee_2)}$
mapping $g=(g_i)$ to $\left(
   \begin{array}{cc}
     g_i & 0 \\
     0 & 1 \\
   \end{array}
 \right)
$. Then we define
$$G_{\bd}=\varinjlim_{\udd(\uee)\in\udim^{-1}(\bd)}G_{\udd(\uee)}.$$ By definition, there are canonical morphisms $f_{\uee}: G_{\uee}\rightarrow G_{\bd}$. The
action of the algebraic group $G_{\udd(\uee)}$ on $\mp^b(A,\uee)$
naturally induces an action of $G_{\bd}$ on $\mp^b(A, \bd).$ Let
$g_{\uee}\in G_{\uee}$ and $x_{\uee'}\in \mp^b(A,\uee').$ Then there
exists $\uee''$ such that $\uee\leq\uee''$ and $\uee'\leq \uee''.$
We define
$$f_{\uee}(g_{\uee}).t_{\uee'}(x_{\uee'})=f_{\uee''}(f_{\uee\uee''}(g_{\uee}).t_{\uee'\uee''}(x_{\uee'})).$$ It is
well-defined. We denote by $\mq\mp^b(A, \uee)$ and
$\mq_1\mp^b(A,\bd)$ the orbit spaces of $\mp^b(A, \uee)$ and
$\mp^b(A, \bd)$ under the actions of $G_{\udd(\uee)}$ and $G_{\bd},$
respectively. We have the follow result.
\begin{Prop}\label{quotientspace} With the above notations, we have
$$
\mq\mp^b(A,\bd)=\mq_1\mp^b(A,\bd)=\varinjlim_{\udd(\uee)\in\underline{\mathrm{dim}}^{-1}(\bd)}\mq\mp^b(A,\uee).
$$
\end{Prop}

The proposition is the direct corollary of the following three
lemmas (see \cite{BD2003} or \cite{JSZ2005}).

\begin{Lemma}\label{lem1} In $\mp^b(A),$ any projective complex can be
uniquely decomposed (up to isomorphism) into the direct sum of a
minimal projective complex and a contractible projective complex.
\end{Lemma}

\begin{Lemma}\label{lem2} Let $f^{\blt}: P^{\blt}\ra Q^{\blt}$ be a morphism
between two minimal projective complex in $\mp^b(A).$ Then
$f^{\blt}$ is a quasi-isomorphism if and only if $f^{\blt}$ is an
isomorphism.
\end{Lemma}

\begin{Lemma}\label{lem3} Let $f^{\blt}: P^{\blt}\ra Q^{\blt}$ be a morphism
in $\mp^b(A)$ and $\uee(P^{\blt})=\uee(Q^{\blt}).$ Then $f^{\blt}$
is a quasi-isomorphism if and only if $f^{\blt}$ is an isomorphism.
\end{Lemma}

\subsection{}
The aim of this subsection is to build the connection between $\mq^b(A, \bd)$ and $\mq\mp^b(A, \bd)$.

\begin{Lemma}  Suppose the following diagram is a pullback of A-module, and
$g_1$ is surjective,

\begin{equation*}
 \begin{CD}
X @>f_1>>Y \\
@VVf_2 V@VVg_1 V \\
Z @>g_2>>W
\end{CD}
    \end{equation*}

then we have the following properties hold:

\begin{enumerate}
  \item $Kerf_1\cong Kerg_2$;

\item $\frac{Y}{Imf_1}\cong \frac{W}{Img_2}$;

\item there exists the exact sequence:
$
 0\longrightarrow X
\longrightarrow Y\oplus Z \longrightarrow W\longrightarrow 0$.

\end{enumerate}
\end{Lemma}
\begin{proof} The first and third statement just follow the
definition of pullback, refer to \cite{ARS1995}. For the second
statement, we use the first statement again to get an
isomorphism:$Kerf_2\cong Kerg_1$. By this, the conclusion follows.
\end{proof}

\nd For any dimension vector sequence $\udd=(\ud_i)_{i\in \bbz},$ we
construct for any complex $M^{\blt}$ in $\mc^b(A,\udd)$ a projective
complex $F^{\blt}$ such that $F^{\blt}$ is quasi-isomorphic to
$M^{\blt}$. Assume $\dim_{\bbc}A=n$ and $\mathrm{gl.dim}. A=m.$

Let $M^{\blt}$ be a complex with the following form:

 \begin{equation}\label{e11}
  \begin{CD}
  0  @>>>M_1 @>\pt_2>>\dots  @>\pt_{r-1}>> M_{r-1} @>\pt_r>> M_{r} @>>> 0
 \end{CD}
 \end{equation}
which $\dim_{\bbc}M_i=d_i$ for any $i\in \bbz.$ Here, if
$\ud_i=(d_i^j)_{j\in Q_0},$ then $d_i=\sum_{j\in Q_0}d_i^j.$

Since $M_r$ has dimension $d_r$, we have the surjective map:
$\pi_r:A^{d_r}\longrightarrow M_r$. Along the differential $\pt_r$
and the above $\pi_r$, we form the pullback $X_{r-1}$:

 \begin{equation}\label{e12}
  \begin{CD}
    \dots @>>>M_{r-2}   @>>> X_{r-1} @>\hat{\pt}_{r}>> A^{d_r} @>>> 0  \\
  &&  @V id VV  @V\hat{\pi}_{r}VV  @V\pi_rVV   \\
  \dots @>>> M_{r-2}@>>> M_{r-1} @>\pt_{r}>> M_r @>>> 0
 \end{CD}
 \end{equation}
Depending on Lemma 2.6, we have the following exact sequence:
$$
 0\longrightarrow X_{r-1}
\longrightarrow A^{d_r}\oplus M_{r-1} \longrightarrow
M_{r}\longrightarrow 0
$$
The dimension of $X_{r-1}$, denoted by $l_{r-1}$, is
$d_{r-1}+d_{r}(n-1)$. Similarly, we have the surjective map:
$\tilde{\pi}_r:A^{l_{r-1}}\longrightarrow X_{r-1}$, and
$\hat{\pi}_r$ is also surjective by the lemma. So we can also form
the pullback $X_{r-2}$ as showed in the following diagram:

\begin{equation}
\xymatrix{\dots \ar[r]&M_{r-3}\ar@{=}[d]\ar[r]& X_{r-2} \ar[r]^{\hat{\pt}_{r-1}} \ar[d]^{\hat{\pi}_{r-1}} & A^{l_{r-1}} \ar[r] \ar[d]^{\tilde{\pi}_{r}} & A^{l_{r}} \ar@{=}[d]\ar[r]& 0\\
 \dots\ar[r]&M_{r-3}\ar[r]\ar@{=}[d]&M_{r-2}\ar[r]\ar@{=}[d]& X_{r-1} \ar[r]^{\hat{\pt}_{r}} \ar[d]^{\hat{\pi}_{r}} & A^{d_{r}} \ar[r] \ar[d]^{\pi_{r}} & 0\\
\dots \ar[r]& M_{r-3}\ar[r]^{\pt_{r-2}}&M_{r-2} \ar[r]^{\pt_{r-1}}
& M_{r-1} \ar[r]^{\pt_r} & M_r \ar[r]& 0}
\end{equation}
Inductively, we get a complex of `almost' free A-module $F^{\blt}$
as follows:
 \begin{equation}\label{e11}
  \begin{CD}
  0 @>>>P@>>>\dots @>>>A^{l_{1}} @>\hat{\pt}_2\tilde{\pi}_2>> \dots @>>>A^{l_r}@>>> 0
 \end{CD}
 \end{equation}
where $l_{i}+l_{i-1}=nl_{i}+d_{i-1}$ for $i=2-m,\cdots, r$,  in
particular, $l_{r}=d_{r}$ and $P$ is a projective module of
dimension $nl_{2-m}-l_{2-m}.$ Every term of this complex is a free
$A$-module except the first term. By the construction, there exists
a projective dimension sequence $\uee$  only depending on the choice of dimension
vector sequence $\udd$ such that $F^{\bullet} \in \mp^b(A,\uee)$. Moreover,
This complex is quasi-isomorphism to $M^{\blt}.$ First,
$$
H_{r}(F^{\blt})=\frac{A^{l_r}}{Im\hat{\pt}_r\tilde{\pi}_r}=\frac{A^{l_r}}{Im\hat{\pt}_r}\cong\frac{M_r}{Im\pt_r}
$$
This follows from that $\tilde{\pi}_r$ is surjective and the above
lemma. In general, for $i<r,$
$$
H_{i}(F^{\blt})=\frac{Ker\hat{\pt}_{i+1}\tilde{\pi}_{i+1}}{Im
\hat{\pt}_i\tilde{\pi}_i}=\frac{Ker\hat{\pt}_{i+1}\tilde{\pi}_{i+1}}{Im\hat{\pt}_i}=\frac{(\hat{\pi}_{i+1}\tilde{\pi}_{i+1})^{-1}(Ker\pt_{i+1})}{(\hat{\pi}_{i+1}\tilde{\pi}_{i+1})^{-1}(Im\pt_i)}=H_{i}(M^{\blt})
.$$ In fact, we can construct a new complex from any place of a
complex and two complexes are quasi-isomorphic to each other as
follows. Given a dimension vector sequence $\udd(i)$ for some $i\in \bbn$, let $M^{\bullet}\in \mc^b(A, \udd(i))$. Then we obtain a complex $X^{\bullet}\in \mc^b(A, \udd(i-1))$ for some dimension vector sequence $\udd(i-1)$ by the commutative diagram
\begin{equation}
\xymatrix{X^{\bullet}: \ar[d]&\cdots\ar[r]\ar@{=}[d]&X_{i-1}\ar[r]\ar[d]&
A^{d_i}\ar[r]\ar[d]&
M_{i+1}\ar[r]\ar@{=}[d]&\cdots\ar@{=}[d]\\
M^{\bullet}:&\cdots\ar[r]& M_{i-1}\ar[r]&M_{i}\ar[r]& M_{i+1}\ar[r]&\cdots}
\end{equation}
where $d_i=\mathrm{dim}_{\bbc}M_i$ and $X_{i-1}$ is the pullback.
 In this way, we obtain a map
$$f_i: \mc^b(A,\udd(i))\rightarrow \mc^b(A,\udd(i-1))$$ such that $f_i(M^{\bullet})$ is quasi-isomorphic to $M^{\bullet}$ for any $M^{\bullet}\in \mc^b(A, \udd(i)).$
\begin{Prop}\label{morphism}
 The above construction induces a map $f_{\udd}: \mc^b(A, \udd)\rightarrow \mp^b(A, \uee)$ for some projective dimension sequence $\uee$ such that there exists a finite stratification $\mc^b(A,\udd)=\bigsqcup_i\mo_i$ such that all $\mo_i$ are constructible and  $f_{\udd}\mid_{\mo_i}$ is a morphism of varieties.
\end{Prop}

The map satisfies the property in Proposition \ref{morphism} is
called the constructible map which has been considered in \cite{Xu}
or \cite{Palu}.

Now we consider the geometrization of Kernel and cokernel of the homomorphisms of modules. First of all, we consider the kernel and cokernel of a linear map.

\begin{Lemma}\label{morphism1} For any $d_1, d_2\in \bbn$, there exist
constructible maps $$\mathcal{K}: Hom_{\bbc}(\bbc^{d_1}, \bbc^{d_2})\rightarrow
\bigcup_{d\leq d_1}Inj(k^{d},k^{d_1})$$ and
$$\mathcal{C}: Hom_{\bbc}(\bbc^{d_1}, \bbc^{d_2})\rightarrow
\bigcup_{d'\leq d_2}Surj(k^{d},k^{d_1})$$ such that for any $f\in Hom_{\bbc}(\bbc^{d_1}, \bbc^{d_2}),$ there exists a long exact sequence
$$
\xymatrix{0\ar[r]&\bbc^{d}\ar[r]^{\mathcal{K}(f)}&\bbc^{d_1}\ar[r]^{f}&\bbc^{d_2}\ar[r]^{\mathcal{C}(f)}&\bbc^{d'}\ar[r]&0}
$$
where $d$ is the rank of $f$ and $d'=d_2+d-d_1.$
\end{Lemma}
It can be easily proved by making a finite partition of $M_{d_2\times d_1}(\bbc).$ As an analogy of this lemma for homorphisms of modules, we have

\begin{Lemma}\label{morphism2}
For any two dimension vectors $\ud_1, \ud_2$, let $(\bbc^{\ud_1}, x_1)$ and $(\bbc^{\ud_2}, x_2)$ be two $A$-modules. Then there exist constructible maps
$$
(\mathcal{K}_1, \mathcal{K}_2): Hom_{A}((\bbc^{\ud_1}, x_1),(\bbc^{\ud_2}, x_2))\rightarrow \bigcup_{\ud; \ud\leq \ud_1}\bbe_{\ud}(A)\times Hom_{\bbc}(\bbc^{\ud}, \bbc^{\ud_1})
$$
and
$$
(\mathcal{C}_1, \mathcal{C}_2): Hom_{A}((\bbc^{\ud_1}, x_1),(\bbc^{\ud_2}, x_2))\rightarrow \bigcup_{\ud'; \ud'\leq \ud_2}\bbe_{\ud'}(A)\times Hom_{\bbc}(\bbc^{\ud_2}, \bbc^{\ud'})
$$
such that for any $f\in Hom_{A}((\bbc^{\ud_1}, x_1),(\bbc^{\ud_2}, x_2))$, there exists a long exact sequence
$$
\xymatrix{0\ar[r]&(\bbc^{\ud}, \mathcal{K}_1(f))\ar[r]^{\mathcal{K}_2(f)}&(\bbc^{d_1}, x_1)\ar[r]^{f}&(\bbc^{d_2}, x_2)\ar[r]^{\mathcal{C}_2(f)}& (\bbc^{d'}, \mathcal{C}_2(f))\ar[r]&0}.
$$
\end{Lemma}
Now we come to prove our proposition.
\\
\begin{proof}
Without loss of generality, we assume $\udd=(\ud_{1}, \cdots,
\ud_r)$ for some $r\in \bbn$. We set $\udd(r):=\udd$. By Lemma
\ref{morphism2}, the above construction of
`almost' free resolution induces a chain of constructible map:
$$
\xymatrix{\mc^b(A, \udd(r))\ar[r]^{f_r}&\mc^b(A,
\udd(r-1))\ar[r]^-{f_{r-1}}&\cdots \ar[r]&\mc^b(A, \udd(-m))}
$$
where $\udd(r-1)=(\ud^{r-1}_{1},\cdots, \ud^{r-1}_r)$ satisfies
$\ud^{r-1}_{i}=\ud_{i}$ for $i<r-1$ and $\ud^{r-1}_{r}=\udim A^{d_r}$ and
$\ud^{r-1}_{r-1}=\udim A^{d_r}+\ud_{r-1}-\ud_r.$
By the construction, $\udd(-m)$ is determined by $\udd$ and $\udim
A$. Let $f_{\udd}=f_{r}\cdots f_{-m}.$ Note that the composition of
constructible maps is  constructible. Then we deduce a constructible
map from $\mc^b(A, \udd(1))$ to $\mc^b(A, \udd(-m))$ satisfying the
image of any complex under $f_{\udd}$ is a almost free projective complex. We complete
the proof of Proposition \ref{morphism}.
\end{proof}

\subsection{} We have the following results.

\bigskip
\nd{\bf Theorem A}  {\em For any $\bd\in K_0(\md^b(A))$, let $\phi_A: \mp^b(A, \bd)\rightarrow \mc^b(A, \bd)$ be the natural embeeding. Then there exists a map $\phi'_A: \mc^b(A, \bd)\rightarrow \mp^b(A, \bd)$ such that for any $\udd\in \udim^{-1}(\bd)$, $\phi'\mid_{\mc^b(A, \udd)}$ is a constructible map and the quotient maps of $\phi_A$ and $\phi'_A$ between $\mq^b(A, \bd)$ and
  $\mq\mp^{b}(A, \bd)$ are inverse to each other.}

\begin{proof}
By Proposition \ref{morphism}, there is a map $\phi'_A$ such that $\phi'_A\mid_{\mc^b(A, \udd)}=f_{\udd}.$
For any $X^{\bullet}\in \mc^b(A, \bd)$, by definition, $X^{\bullet}$ is quasi-isomorphic to $\phi_A\phi'_A(X^{\bullet}).$  This proves the theorem.
\end{proof}

\bigskip
Let us  recall some results  on Morita theory of derived categories (\cite{Rickard1989} and
\cite{Rickard1991}). First,
a tilting complex $T$ over $A$ is an object in $\mk^b(P)$ which
satisfies the following conditions:

1. For any $i\neq 0$, $\mathrm{Hom}_{\md^b(A)}(T,T[i])=0$;

2. The category $\mathrm{add}(T)$ generates $\mk^b(P)$ as a
triangulated category.

\nd For convenience, we shall consider the right A-module
temporarily and $A$-$B$-bimodule means module for
$A^{op}\bigotimes_{\bbc}B$ in the following part of this section.
J. Rickard proved the following result.

\nd{\bf Rickard's Theorem}  {\em For any two finite dimensional basic $\bbc$-algebras $A$ and $B,$ they are derived equivalent
 if and only if there is a tilting complex $T$ over $A$ such that the endomorphism ring
$\mathrm{End}_{\md^b(A)}(T)^{op}$ is isomorphic to $B.$}

\bigskip
Moreover, the method of Rickard implies the further results. If
$A$ and $B$ are derived equivalent by a functor $F,$ then the
functor $F$ induces  a derived equivalence between
$\md^b(A^{op}\otimes_{\bbc}A)$ and $\md^b(B^{op}\otimes_{\bbc}A)$.
The image of $A$ as $A^{op}\otimes_{\bbc}A$ module under this
functor is a complex,  say $\Delta$, in
$\md^b(B^{op}\otimes_{\bbc}A)$. The functor $F$ induces also a
derived equivalence between $\md^b(B^{op}\otimes_{\bbc}B)$ and
$\md^b(A^{op}\otimes_{\bbc}B)$. The image of $B$ as
$B^{op}\otimes_{\bbc}B$ module under this functor is a complex,
say $\Theta$, in $\md^b(A^{op}\otimes_{\bbc}B).$ Then
$$
-\otimes_{B}^{L}\Delta: \md^b(B)\longrightarrow \md^b(A)
$$
is an equivalence of triangulated categories with two
quasi-inverses
$
\mathrm{RHom}_{\md^b(A)}(X,-): \md^b(A)\longrightarrow \md^b(B)
$
and
$
-\otimes_{B}^{L}\Theta: \md^b(A)\longrightarrow \md^b(B).
$

\bigskip

\nd{\bf Theorem B} {\em Let $F$: $\md^b(A)\rightarrow \md^b(B)$ be
a functor of derived equivalence with the induced isomorphism $F_0: K_0(\md^b(A))\rightarrow K_0(\md^b(B))$. Then for any $\bd_A\in K_0(\md^b(A))$, there exist maps $\mathcal{F}_1: \mc^b(A, \bd_A)\rightarrow \mc^b(B, F_0(\bd_A))$ and $\mathcal{F}_2: \mc^b(B, F_0(\bd_A))\rightarrow \mc^b(A, \bd_A)$ such that for any $\udd_A\in \udim^{-1}(\bd_A)$ and $\udd_B\in \udim^{-1}(F_0(\bd_A))$, $\mathcal{F}_1\mid_{\mc^b(A, \udd_A)}$ and $\mathcal{F}_2\mid_{\mc^b(B, \udd_B)}$ are constructible maps and the quotient maps of $\mathcal{F}_1$ and $\mathcal{F}_2$ are inverse to each other.}

\begin{proof}  We consider $\Delta$ as a $B^{op}\bigotimes_{\bbc}A$-projective
complex:

\begin{equation}\label{e11}
  \begin{CD}
  \dots  @>>> \triangle_{j} @>\pt^{\Delta}_{j}>> \Delta_{j+1} @>>> \dots
 \end{CD}
 \end{equation}
where $\Delta_{j}$ is $B^{op}\bigotimes_{\bbc}A$-projective module. For any $X^{\bullet}\in \mc^b(B, \udd_B),$  let $f_{\udd_B}(X^{\bullet})$ is the corresponding
$B$-projective complex with the following form:

\begin{equation}\label{e11}
  \begin{CD}
  \dots  @>>> P^{i} @>\pt^{P}_{i}>> P^{i+1} @>>> \dots
 \end{CD}
 \end{equation}
Here, $f_{\udd_B}$ is the constructible map from $\mc^b(B, \udd_B)$ to $\mp^b(B, \uee_1)$ for some projective dimension sequence $\uee_1$. Define
the following complex (see \cite{Rickard1991}):

\begin{equation}\label{e11}
  \begin{CD}
  \dots  @>>> \prod_{i+j=n}P^{i}\bigotimes_{B}\Delta_{j} @>d_{n}>> \prod_{i+j=n+1}P^{i}\bigotimes_{B}\Delta_{j} @>>> \dots
 \end{CD}
 \end{equation}
where $d_{n}=\pt^{P}_{i}\bigotimes
id+(-1)^{i}id\bigotimes\pt^{\Delta}_{i}$ and
$\prod_{i+j=n}P^{i}\bigotimes_{B}\Delta_{j}$ is $A$-projective
module. We set $\uee_{2}$ to be its $A$-projective dimension sequence.  This defines  a map $f_{\uee_1\uee_2}: \mp^b(B,\uee_{1})\longrightarrow\mp^b(A,\uee_{2})$.
By the above construction, it is clearly a constructible map. In the same way, we can define the constructible map
$$
g_{\uee_2\uee_3}:\mp^b(A,\uee_{2})\longrightarrow\mp^b(B,\uee_{3})
$$
by $g(P^{\blt})=P^{\blt}\bigotimes_{A}\Theta,$  which is induced
by $G=-\bigotimes_{A}\Theta$. The constructible maps $\phi_Bg_{\uee_2\uee_3}$ and $f_{\uee_1\uee_2}f_{\udd_B}$ induce the maps $\mathcal{F}_1$ and $\mathcal{F}_2.$ Because $F$ and $G$ are
quasi-inverse to each other, $\phi_Bg_{\uee_2\uee_3}f_{\uee_1\uee_2}(f_{\udd_B}(X^{\bullet}))$ is isomorphic to $X^{\bullet}$ in $\md^b(B).$
This proves the theorem.
\end{proof}

\section{Constructible functions on topological spaces attached to derived categories}

\subsection{Degenerations in derived categories}
We rewrite the definition of the degeneration in  \cite{JSZ2005} for
our situation.  For any $X$ and $Y$ in $\mk^{b}(\mp^b(A))$ we denote
by $X\leqslant_{\triangle}Y$ if there is a distinguished triangle
$$
\xymatrix{Y \ar[r]  & X\bigoplus{Z} \ar[r] & Z \ar[r]& Y[1]}
$$
for some $Z\in \mk^{b}(\mp(A))$. On the topological side, we denote
by $X\leqslant_{top}Y$ if $ Y\in \overline{\mo(X)}$ in
$\mp^b(A,\bd)$ where $\overline{\mo(X)}$ is the closure of the orbit
$\mo(X)$ of $X$ in $\mp^b(A, \bd)$ under the action of $G_{\bd}.$ In
order to avoid confusion, we use the following different notation:
$X\leqslant_{top}^{*}Y$ if and only if $Y\in
\overline{G_{\uee}X}=\overline{\mo_{\uee}(X)}$ in $ \mp^b(A,\uee)$
for a fixed projective dimension sequence $\uee$ where
$\mo_{\uee}(X)$ is the orbit of $X$ in $\mp^b(A, \uee)$ under the
action of $G_{\uee}$ and $\overline{\mo(X)}$ is the closure of
$\mo(X)$ in $\mp^b(A,\uee).$ Then Theorem {\bf 1} and Theorem {\bf
2} in \cite{JSZ2005} implies the following result.

\begin{theorem}
For any $X$ and $Y$ in $\mk^{b}(\mp(A))$, then
$X\leqslant_{\triangle}Y$ if and only if $X\leqslant_{top}Y$.
\end{theorem}

For any projective dimension sequence $\uee'$ such that $\uee\leq
\uee'$, we set $X_{\uee'}:=t_{\uee\uee'}(X)$. Then we have
$$
\mo(X)=\varinjlim_{\udd(\uee')\in\underline{\mathrm{dim}}^{-1}(\bd)}\mo_{\uee'}(X_{\uee'}),
\quad\quad\overline{\mo(X)}=\varinjlim_{\udd(\uee')\in\underline{\mathrm{dim}}^{-1}(\bd)}\overline{\mo_{\uee'}(X_{\uee'})}.
$$

\begin{Prop}
Let $X\in \mp^b(A, \uee)$ where $\uee\in
\underline{\mathrm{dim}}^{-1}(\bd)$. Then $\mo(X)$ is a locally
closed subset of $\mp^b(A,\bd)$, i.e., the intersection of a closed
subset with an open subset.
\end{Prop}
\begin{proof} For any $\uee'\geq\uee,$ consider the morphism $t_{\uee\uee'}: \mp^b(A, \uee)\rightarrow \mp^b(A,\uee')$
sending $M^{\bullet}$ to $M^{\bullet}\oplus P^{\bullet}(\uee,\uee')$
where $P^{\bullet}(\uee,\uee')$ is the direct sum of the complexes
with only two nonzero term formed as $P\xrightarrow{f} P.$ For any
$\lambda\in \bbc^*,$ we denote by $P_{\lambda}^{\bullet}(\uee,
\uee')$ the complex isomorphic to $P^{\bullet}(\uee, \uee')$
obtained by substituting $P\xrightarrow{\lambda f} P$ for any direct
summand $P\xrightarrow{f} P$ of $P^{\bullet}(\uee, \uee')$.  The
orbit $\mo_{\uee}(X)$ of $X$ in $\mp^b(A, \uee)$ is locally closed.
By definition, $\mo_{\uee}(X)=\overline{\mo_{\uee}(X)}\bigcap
U_{\uee}$ for some open subset $U_{\uee}$ in $\mp^b(A, \uee)$. It is
clear that the set $\mathcal{S}(\uee, \uee'):=\{X^{\blt}\oplus
P_{\lambda}^{\bullet}(\uee, \uee')\mid X^{\blt}\in U_{\uee},
\lambda\in \bbc^* \}$ is an open subset of $\mp^b(A, \uee')$. Set
$U_{\uee'}:=\bigsqcup_{g\in G_{\uee'}}g.\mathcal{S}(\uee, \uee').$
It is an open subset of $\mp^b(A, \uee')$. Then
$\mo_{\uee'}(X_{\uee'})=\overline{\mo_{\uee'}(X_{\uee'})}\bigcap
U_{\uee'}$.   Then by definition, $\mo(X)$ is equal to
$$
\varinjlim_{\udd(\uee')\in\underline{\mathrm{dim}}^{-1}(\bd)}\mo_{\uee'}(X_{\uee'})=\varinjlim_{\udd(\uee')\in\underline{\mathrm{dim}}^{-1}(\bd)}\overline{\mo_{\uee'}(X_{\uee'})}\bigcap
U_{\uee'}=\overline{\mo(X)}\bigcap
\varinjlim_{\udd(\uee')\in\underline{\mathrm{dim}}^{-1}(\bd)}U_{\uee'}.
$$
By the definition of the fine topology on $\mp^b(A, \bd)$,
$\varinjlim_{\udd(\uee')\in\underline{\mathrm{dim}}^{-1}(\bd)}U_{\uee'}$
is an open subset of $\mp^b(A, \bd).$ We finish the proof.
\end{proof}
 We recall the following
property which is proved in \cite{JSZ2005}: A complex $P(x)$
corresponding to $x\in \mp^b(A,\uee)$ is partial tilting complex if
$\hom_{\mk^{b}(A)}(P(x),P(x)[1])=0.$ Then
$\mo_{\uee}(x)=\{y\in\mp^b(A,\uee)|y\sim x\}$ is open in
$\mp^b(A,\uee).$ Furthermore, the orbit  $\mo(x)$ in $\mp^b(A,\bd)$
is also open if the corresponding complex $P(x)$ is a partial
tilting complex.

\begin{Prop}
Any point in $\mq\mp^b(A,\bd)$ is not closed, that is, any orbit is
not closed in $\mp^b(A,\bd).$
\end{Prop}
\begin{proof} Take  any projective $A$-module $P$ and
$t\in\bbc,$ we define the complex $C(t)=(P_i,\pt_i)_{i\in\bbz}$ by
$P_0=P_1=P$ and other $P_i=0\ \mbox{for}\ i\neq 0,1;$ $\pt_0=t$ and
other $\pt_i=0\ \mbox{for}\ i\neq0.$

If $t\neq 0$, then $C(t)$ is contractible but $C(0)$ is not
quasi-isomorphic to zero. Moreover, $C(1)\leqslant_{top} {C(0)}$.
 Similarly, For any $x\in \mp^b(A,\bd)$, $P(x)$ is the corresponding
 complex, $P(x)\bigoplus{C(1)}\leqslant_{top} {P(x)\bigoplus{C(0)}}$,
i.e.the point corresponding to $P(x)\bigoplus{C(0)}$ is in the
closure of the point corresponding to $P(x)\bigoplus{C(1)}$(it is
quasi isomorphic to $P(x)$), but not in its orbit. This shows the
orbit of $x$ is not closed. \end{proof}

By this proposition, we can construct an infinite sequence of
non-trivial  degenerations:
$$
P(x)\leqslant_{top}P(x)\bigoplus{C(0)}\leqslant_{top}P(x)\bigoplus{C(0)}\bigoplus{C(0)}
\leqslant\cdots
$$

\subsection{The na\"\i ve Euler characteristic}
Let $X$ be a algebraic variety over $\bbc$.  We denote by $M(X)$ the
set of all constructible functions on algebraic variety $X$ with
values in $\bbc$. The set $M(X)$ is naturally a $\bbc$-linear space.
Let $G$ be an algebraic group acting on $X.$ Then we denote by
$M_{G}(X)$ the subspace of $M(X)$ consisting of all $G$-invariant
functions.

Let $\chi$ denote Euler characteristic in compactly-supported
cohomology. Let $X$ be an algebraic variety and $\mo$ a
constructible subset as the disjoint union of finitely many
locally closed subsets $X_i$ for $i=1,\cdots,m.$ Define
$\chi(\mo)=\sum_{i=1}^m\chi(X_i).$ We note that it is
well-defined. We will use the following properties:
\begin{Prop}[\cite{Dimca},\cite{Riedtmann1994} and \cite{Joyce}]\label{Euler} Let $X,Y$ be algebraic varieties over $\mathbb{C}.$
Then
\begin{enumerate}
    \item  If the algebraic variety $X$ is the disjoint union of
finitely many constructible sets $X_1,\cdots,X_r$, then
$$\chi(X)=\sum_{i=1}^{r}{\chi(X_i)}.$$
    \item  If $\varphi:X\longrightarrow Y$ is a morphism
with the property that all fibers have the same Euler
characteristic $\chi$, then $\chi(X)=\chi\cdot \chi(Y).$ In
particular, if $\varphi$ is a locally trivial fibration in the
analytic topology with fibre $F,$ then $\chi(X)=\chi(F)\cdot
\chi(Y).$
    \item $\chi(\bbc^n)=1$ and $\chi(\mathbb{P}^n)=n+1$ for all $n\geq
    0.$
\end{enumerate}
\end{Prop}
We recall the definition of {\it pushforward} functor from the
category of algebraic varieties over $\mathbb{C}$ to the category of
$\mathbb{Q}$-vector spaces.

\nd Let $\phi: X\rightarrow Y$ be a morphism of varieties. For
$f\in M(X)$ and $y\in Y,$ define
$$
\phi_{*}(f)(y)=\sum_{c\in\bbq}c\chi(f^{-1}(c)\cap \phi^{-1}(y)).
$$
\begin{theorem}[\cite{Dimca},\cite{Joyce}]\label{Joyce}
Let $X,Y$ and $Z$ be algebraic varieties over $\mathbb{C},$ $\phi:
X\rightarrow Y$ and $\psi: Y\rightarrow Z$ be morphisms of
varieties, and $f\in M(X).$ Then $\phi_{*}(f)$ is constructible,
$\phi_{*}: M(X)\rightarrow M(Y)$ is a $\mathbb{Q}$-linear map and
$(\psi\circ \phi)_{*}=(\psi)_{*}\circ (\phi)_{*}$ as
$\mathbb{Q}$-linear maps from $M(X)$ to $M(Z).$
\end{theorem}
In order to deal with orbit spaces, we  need to consider
geometric quotients.
\begin{Definition}
Let $G$ be an algebraic group acting on a variety $X$ and
$\phi:X\rightarrow Y$ be a $G$-invariant morphism, i.e. a morphism
constant on orbits. The pair $(Y,\phi)$ is called a geometric
quotient if $\phi$ is open and for any open subset $U$ of $Y$, the
associated comorphism identifies the ring $\mo_{Y}(U)$ of regular
functions on $U$ with the ring $\mo_{X}(\phi^{-1}(U))^{G}$ of
$G$-invariant regular functions
 on $\phi^{-1}(U)$.
\end{Definition}

The following result due to Rosenlicht \cite{Rosenlicht1963} is
essential to us.

\begin{Lemma}
Let $X$ be a $G$-variety, then there exists a  open and dense
$G$-stable subset which has a geometric $G$-quotient.
\end{Lemma}

By this Lemma, we can construct a finite stratification over $X.$
Let $U_1$ be an open and dense $G$-stable subset of $X$ as in Lemma
3.8. Then we have a geometric $\phi_{U_1}: U_1\rightarrow Y_1$.
Since $\dim_{\mathbb{C}}(X-U_1)<\dim_{\mathbb{C}}X,$ we can use the
above lemma again, there exists a dense open $G$-stable subset $U_2$
of $X-U_1$ which has a geometric $G$-quotient $\phi_{U_2}:
U_2\rightarrow Y_2$. Inductively, we get a finite stratification
$X=\cup_{i=1}^{l}U_{i}$ where $U_{i}$ is a $G$-invariant locally
closed subset and has a geometric quotient $\phi_{U_i}:
U_i\rightarrow Y_i$ for $i=1, \cdots, l$ with $l\leq
\dim_{\mathbb{C}}X.$  Depending on this stratification, we define
the na\"\i ve Euler characteristic of the orbit space of $X$ under
the group action of $G$ as follows:
$$
\chi^{\mathrm{na}}(X/G):=\chi^{\mathrm{na}}([X/G](\bbc))=\chi(\coprod_{i=1}^l
Y_i)=\sum_{i=1}^{l}\chi(\phi_{U_{i}}(U_{i}))
$$
where $[X/G]$ is the quotient stack with the set of $\bbc$-points
$[X/G](\bbc)$ and $[X/G](\bbc)$ is pseudoisomorphic to
$\coprod_{i=1}^l Y_i$ (see \cite[Section 4]{Joyce} for the
definition of pseudomorphism). For simplicity, we substitute $\chi$
for $\chi^{\mathrm{na}}$ in the following. Indeed, when $G=id,$ then
the na\"\i ve Euler characteristic is just the Euler characteristic.
So the overlapping notation should not cause any confusion. It is
well-defined by the following observation.

Let $[X/G](\bbc)$ be the set of $\bbc$-points. A subset $C\subseteq
[X/G](\bbc)$ is constructible if $C=\bigcup_{i\in
I}\mathcal{F}_i(\bbc)$ where $\{\mathcal{F}_i\mid i\in I\}$ is a
finite set of algebraic $\bbc-$ substacks of finite type of $[X/G]$
(see \cite[Definition 4.1]{Joyce}).
\begin{Lemma}\label{Eulersum}
Let $X$ be a $G$-variety. If $[X/G](\bbc)$ is the disjoint union of
finitely many constructible subsets $Z_{1},\cdots,Z_{r}$, then
$\chi(X/G)=\sum_{i=1}^{r}\chi(Z_i)$.
\end{Lemma}
\begin{proof}  By the above construction,
there is a pseudoisomorphism $\Phi$ between $[X/G](\bbc)$ and
$Y:=\coprod_{i=1}^l Y_i$. For $i=1, \cdots, r$,  $\Phi(Z_i)$ is a
constructible subset of $Y$ and $Y=\bigsqcup_{i=1}^r\Phi(Z_i)$. Then
by the definition of the naive Euler characteristic and Proposition
\ref{Euler}, we deduce $\chi(X/G)=\sum_{i=1}^{r}\chi(Z_i).$
\end{proof}

Let $X$ and $Y$ be two complex algebraic varieties under the actions
of the algebraic groups $G$ and $H$, respectively. Let $\phi:
[X/G]\rightarrow [Y/H]$ be a $1$-morphism. Then it induces a natural
pseudomorphism $\phi_*: [X/G](\bbc)\rightarrow [Y/H](\bbc)$
\cite[Proposition 4.6]{Joyce}. In the same way as the proof of Lemma
\ref{Eulersum}, we obtain the following result.
\begin{Lemma}\label{Eulermulti}If $\phi_*$ is surjective and all fibres of
$\phi_*$ have the same naive Euler characteristic $\chi,$ then we
have
$$
\chi(X/G)=\chi(Y/H)\cdot \chi.
$$
\end{Lemma}
 We introduce the following notation. Let $X$ be a variety under
 the action of an algebraic group $G$. Let $f$ be a $G$-invariant
 constructible function over $X.$ We define
 $$
\int_{x\in [X/G](\bbc)}f(x):=\sum_{c\in \bbc}\chi(f^{-1}(c)/G)c.
 $$

\subsection{}
Let $\bd_1, \bd_2$ be two dimension vector in $K_0$. For any two
subsets $\mo_1$ and $\mo_2$ of $\mp^b(A,{\bd}_1)$ and
$\mp^b(A,{\bd}_2)$ respectively, we define the subset $\mo_1*\mo_2$
of $\mp^b(A,{\bd}_1+{\bd}_2)$ to be the set of
$z\in\mp^b(A,{\bd}_1+{\bd}_2)$ such that there exists a triangle
$$
P(y)\longrightarrow {P(z)}\longrightarrow {P(x)}\longrightarrow
{P(x)[1]}
$$ in $\md^b(A)$ where $x\in\mo_1$ and $y\in\mo_2.$ By the
octahedral axiom, we have $(\mo_1*\mo_2)*\mo_3=\mo_1*(\mo_2*\mo_3).$
Inductively, we can define $\mo_1*\mo_2*\cdots*\mo_s$ for all $s>1.$
A subset $\mo$ of $\mp^b(A,\bd)$ is called a support-bounded
constructible subset if there exists a constructible subset
$\mo_{\uee}$ of $\mp^b(A, \uee)$ for some $\udd(\uee)\in
\udim^{-1}(\bd)$ such that $\mo=G_{\bd}.t_{\uee}(\mo_{\uee})$ where
$t_{\uee}: \mp^b(A, \uee)\rightarrow \mp^b(A, \bd)$ is canonical.
Let $\udd(\uee')\in \udim^{-1}(\bd)$ and $\uee\leq \uee'.$ We set
$\mo_{\uee'}:=G_{\uee'}.t_{\uee\uee'}(\mo_{\uee}).$ Then we have
$$
\mo=\varinjlim_{\udd(\uee')\in\underline{\mathrm{dim}}^{-1}(\bd)}\mo_{\uee'}.
$$

 The following proposition shows the property of
support-bounded is invariant under derived equivalence.
\begin{Prop}
If two algebras $A$ and $B$ are derived equivalent, then this
equivalent functor $F$ induces the map $f$ as in Theorem B sending
constructible sets of support bounded in $\mp^b(A)$ to
constructible sets of support bounded in $\mp^b(B)$.
\end{Prop}
\begin{proof} By Theorem $B$, the equivalent functor
$F:D^{b}(A)\rightarrow D^{b}(B)$ induces the following commutative
diagram:
\begin{equation}\label{e12}
  \begin{CD}
   \mp^{b}(A,\bd) @>f >> \mp^{b}(B,\bd')   \\
    @V\pi_{A} VV  @V\pi_{B} VV \\
   \mq\mp^b(A,\bd) @>\bar{f} >> \mq\mp^b(B,\bd')
  \end{CD}
 \end{equation}
 where $f$ is a morphism of varieties and $\bar{f}$ is a
 homeomorphism of quotient spaces. Let $\mo$ be a support-bounded constructible subset of $\mp^b(A,\bd).$
  Then $\mo=G_{\bd}.\mo_{\uee}$  for some projective dimension sequence $\uee$ by
definition. We have $f(\mo_\uee)$ is a constructible subset of $
\mp^{b}(B,\uee')$ for some projective dimension sequence $\uee'$.
Since $F$ is a derived equivalence,  $f(\mo)$ is
$G_{\bd'}$-invariant. On the other hand, due to the commutativity of
the diagram,
$$
f(\mo)=\pi_{B}^{-1}\overline{f}\pi_{A}(\mo)=\pi_{B}^{-1}\overline{f}\pi_{A}(\mo_\uee)=G_{\bd'}.f(\mo_{\uee}).
$$
\end{proof}

\begin{Definition}\label{equivalence class-1}Let $\bd$ be a dimension vector in $K_0$ ($\mk_0(\md^b(A))$) and
$\mo$ be a support-bounded constructible subset of $\mp^b(A,\bd).$
Then we denote by $\hat{1}_{\mo}$ the $\bbc$-value function over
$\mp^b(A, \bd)$ given by taking values $1$ on  each point in $\mo$
and $0$ otherwise. A function $\hat{f}$ on $\mp^b(A,\bd)$ is called
a $G_{\bd}$-invariant constructible function
 if $\hat{f}$ can be written as a sum of finite terms $\sum_{i}m_i \hat{1}_{\mo_i}$
 where $m_i\in \bbc$ and $\mo_i$ is a support-bounded constructible subset of $\mp^b(A,\bd)$.
The set of $G_{\bd}$-invariant constructible functions over
$\mp^b(A,\bd)$ is denoted by $M(\mp^b(A,\bd)).$
\end{Definition}
 Let $\uee\in \udim^{-1}(\bd)$ and $M(\mp^b(A,
\uee))$ be the $\bbc$-space of $G_{\uee}$-invariant constructible
functions over $\mp^b(A, \uee).$ For any $\uee' \in \udim^{-1}(\bd)$
and $\uee'\geq \uee,$ there is a natural linear map $r_{\uee'\uee}:
M(\mp^b(A, \uee'))\rightarrow M(\mp^b(A, \uee))$ sending a function
$\hat{f}$ to $t_{\uee\uee'}\circ \hat{f}.$ Then we obtain an inverse
system.
\begin{Prop}
With the notations in Definition \ref{equivalence class-1}, we have
$$
M(\mp^b(A,\bd))=\varprojlim_{\udd(\uee)\in\underline{\mathrm{dim}}^{-1}(\bd)}M(\mp^b(A,
\uee)).
$$
\end{Prop}

\subsection{Convolution}

Let $\bd', \bd''$be two dimension vectors in $K_0$ and $\uee',
\uee''$ be two projective dimension vectors satisfying
$\udd(\uee')\in \udim^{-1}(\bd'), \udd(\uee'')\in
\udim^{-1}(\bd'')$. The following two lemmas is easy corollaries of
4.1.13 in \cite{GM1996}.
\begin{Lemma}\label{lem3.1}
If there is a distinguished triangle: $ \xymatrix{X^{\blt}
\ar[r]^{\underline{f}} & L^{\blt} \ar[r]^{\underline{g}}  & Y^{\blt}
\ar[r]^{\underline{h}} & X^{\blt}[1]} $ for $X^{\blt}\in
\mp^{b}(A,\uee''),$ $Y^{\blt}\in \mp^{b}(A,\uee')$, and $L^{\blt}\in
\mp^{b}(A,\uee'+\uee'')$, then it induces an exact sequence of chain
complex:
$$
\xymatrix{0 \ar[r]& X^{\blt} \ar[r] & L^{\blt} \ar[r]  & Y^{\blt}
\ar[r] & 0}
$$
\end{Lemma}

\begin{Lemma}\label{lem3.2}
Let $X^{\blt}\in \mp^{b}(A,\uee^{''}),Y^{\blt}\in
\mp^{b}(A,\uee^{'})$ and the exact sequence:
 $$
\xymatrix{0 \ar[r]& X_{i}
\ar[r]^{l_{i}} & L_{i} \ar[r]^{\pi_{i}}  & Y_{i} \ar[r] & 0}
$$
hold for all   $i$-th terms of  $X^{\blt}$, $Y^{\blt}$ and
$L^{\blt}$ respectively, where $l_i$ are the canonical injections
and $\pi_i$ the canonical projections. Define the set $\{\pt^L\}$ to be the set of sequences
$(\pt_i^L)_i$ such that $L^{\blt}=(X_i\oplus Y_i, \pt_i^L)_i$
becomes a complex and  induces the exact sequence of complexes

$$
\xymatrix{0 \ar[r]^{}& X^{\blt} \ar[r]^{} & L^{\blt} \ar[r]^{} &
Y^{\blt} \ar[r] & 0}.
$$Then we have the canonical
isomorphism $\{\pt^{L}\}\cong
\hom_{\mc^b(A)}(Y^{\blt}[-1],X^{\blt}).$
\end{Lemma}

Let $\mo_1\subset\mp^b(A,\uee'')$ and $\mo_2\subset\mp^b(A,\uee')$
be $G_{\uee''}$- and $G_{\uee'}$-invariant constructible set,
respectively. Put
$$
L(\mo_1,\mo_2)=\{Cone(h)\mid h\in
\hom_{\mk^b(\mp(A))}(\mo_1[-1],\mo_2)\}.
$$
For $L\in\mp^b(A,\uee'+\uee'')$, we set
$$
W_{\uee'\uee''}(\mo_1,\mo_2;L)=\{(f,g,h)\mid
\xymatrix{Y \ar[r]^{f}& L \ar[r]^{g} & X \ar[r]^{h} & Y[1] }\\
\mbox{is a distinguished triangle }, $$
$$
\mbox{with } X\in \mo_1,Y\in\mo_2\}
$$
$$
W_{\uee'\uee''}(\mo_1,\mo_2)=\{(L, (f,g,h))\mid
L\in\mp^b(A,\uee'+\uee''), (f,g,h)\in W_{\uee'\uee''}(\mo_1,\mo_2;L)
\}$$ and
$$
\hom_{\mk^b(\mp(A))}(\mo_1[-1],\mo_2)=\{(Y,X,h)\mid X\in\mo_1,\
Y\in\mo_2, h\in\hom_{\mk^b(\mp(A))}(X[-1],Y) \}.
$$
 We introduce the action of $G_{\uee''}\times
G_{\uee'}$ on $W_{\uee'\uee''}(\mo_1,\mo_2)$
and $\hom_{\mk^b(\mp(A))}(\mo_1[-1],\mo_2)$ as follows:\\
For $(a,c)\in G_{\uee''}\times G_{\uee'},$
$$(a,c)\circ (L,(f,g,h))=(L,(fc^{-1},ag,(c[1])ha)),$$
$$(a,c)\circ (h)=ch(a[-1])^{-1}.$$

The action of $G_{\uee''}\times G_{\uee'}$ on
$W_{\uee'\uee''}(\mo_1,\mo_2)$ naturally induces the action on
$W_{\uee'\uee''}(\mo_1,\mo_2;L)$. We consider the na\"\i ve Euler
characteristic of quotient space
$W_{\uee'\uee''}(\mo_1,\mo_2;L)/{G_{\uee''}\times G_{\uee'}}$ and
have the following result.
\begin{Prop}\label{constructible}
Let $\mo_{1},\mo_{2}$  be as above. Then the set
$$\{\chi(W_{\uee'\uee''}(\mo_1,\mo_2;L)/{G_{\uee''}\times G_{\uee'}})\mid L\in
\mp^b(A,\uee'+\uee'')\}$$ is a finite set.
\end{Prop}

\begin{proof}  Consider the quotient stack $[W_{\uee'\uee''}(\mo_1,\mo_2)/G_{\uee''}\times G_{\uee'}]$
 of $W_{\uee'\uee''}(\mo_1,\mo_2)$ under the action of $G_{\uee''}\times
 G_{\uee'}.$ As in Section 3.2, there exists a pseudoisomorphism $\Phi$
 between $[W_{\uee'\uee''}(\mo_1,\mo_2)/G_{\uee''}\times G_{\uee'}](\bbc)$ and some
 variety $Y.$ The natural projection
$\pi: W_{\uee'\uee''}(\mo_1,\mo_2)\rightarrow \mp^b(A,\uee'+\uee'')$
and $\Phi^{-1}$ combines to give a pseudomorphism $\Psi:
Y\rightarrow \mp^b(A,\uee'+\uee'').$ It is a morphism of varieties
\cite[Section 3.4]{Joyce}. By Theorem \ref{Joyce}, $\pi_{*}(1_{Y})$
is constructible. By the defintion of pushforward functor, it means
that the set
$$\{\chi(\Psi^{-1}(L))\mid L\in
\mp^b(A,\uee'+\uee'')\}$$ is a finite set. Therefore the set
$\{\chi(W_{\uee'\uee''}(\mo_1,\mo_2;L)/{G_{\uee''}\times
G_{\uee'}})\mid L\in \mp^b(A,\uee'+\uee'')\}$ is finite.
\end{proof}

\begin{Definition}\label{convolution-1}
Let $\mo_1\subset\mp^b(A,\uee'')$ and $\mo_2\subset\mp^b(A,\uee')$
be $G_{\uee''}$- and $G_{\uee'}$-invariant constructible sets,
respectively. The convolution multiplication
$1_{\mo_{1}}*1_{\mo_{2}}\in M(\mp^{b} (A,\uee{'}+\uee{''}))$ is
defined as follows:
$$
1_{\mo_{1}}*1_{\mo_{2}}(L)=\chi(W_{\uee'\uee''}(\mo_1,\mo_2;L)/{G_{\uee''}\times
G_{\uee'}})
$$
for $L\in\mp^b(A,\uee'+\uee'').$ We set $V(\mo_1,
\mo_2):=[W_{\uee'\uee''}(\mo_1,\mo_2)/G_{\uee''}\times G_{\uee'}]$
and $V_{\uee'\uee''}(\mo_1,
\mo_2;L):=[W_{\uee'\uee''}(\mo_1,\mo_2;L)/G_{\uee''}\times
G_{\uee'}].$ Write
$F_{\mo_1,\mo_2}^{L}=\chi(V_{\uee'\uee''}(\mo_1,\mo_2;L))$ .
\end{Definition}
Obviously $1_{\mo_{1}}*1_{\mo_{2}}$ is again
$G_{\uee'+\uee''}$-invariant. In this way, the proof of Proposition
\ref{constructible} implies
\begin{Cor}\label{constructible-well-defined}
If $f\in M(\mp^b(A,\uee''))$ and $g\in M(\mp^b(A,\uee')),$ then the
convolution $f*g\in M(\mp^b(A,\uee'+\uee''))$ is well-defined.
\end{Cor}

The above discussion can be extended to $\mp^b(A, \bd).$ Assume that
$\uee'\in \udim^{-1}(\bd')$ and $\uee''\in \udim^{-1}(\bd'')$. Put
$\mo=G_{\bd''}.t_{\uee''}(\mo_1)$ and
$\mo'=G_{\bd'}.t_{\uee'}(\mo_2).$ Let $\uee'_0\geq \uee'$ and
$\uee''_0\geq \uee''$ be in $\udim^{-1}(\bd')$ and
$\udim^{-1}(\bd'')$, respectively. Then $t_{\uee'\uee'_0}$ and
$t_{\uee''\uee''_0}$ naturally induces a map $t_{\uee'\uee'_0,
\uee''\uee''_0}: W_{\uee'\uee''}(\mo_1,\mo_2)\rightarrow
W_{\uee'_0\uee''_0}(\mo_{1\uee''_0}, \mo_{2\uee'_0})$ where
$\mo_{1\uee'_0}=G_{\uee'_0}.t_{\uee'\uee'_0}(\mo_1)$ and
$\mo_{2\uee''_0}=G_{\uee''_0}.t_{\uee''\uee''_0}(\mo_2)$. Then we
obtain a direct system and set
$$
W(\mo, \mo')=\varinjlim W_{\uee'_0\uee''_0}(\mo_{1\uee''_0},
\mo_{2\uee'_0}).
$$
For any $L\in \mp^b(A, \bd'+\bd''),$ we set
$L_{\uee}:=t_{\uee}^{-1}(L)\in \mp^b(A, \uee)$ (perhaps it does not
exist!). Then we define
$$
W(\mo, \mo'; L)=\varinjlim W_{\uee'_0\uee''_0}(\mo_{1\uee''_0},
\mo_{2\uee'_0}; L_{\uee'_0+\uee''_0}).
$$
There are natural actions of $G_{\bd''}\times G_{\bd'}$ on $W(\mo,
\mo'; L)$ and $W(\mo, \mo').$ The orbit spaces are denoted by
$V(\mo, \mo'; L)$ and $V(\mo, \mo'),$ respectively. We note that
$V(\mo,\mo'; L)=V_{\uee'\uee''}(\mo_1,\mo_2; L_{\uee'+\uee''}).$

\begin{Definition}\label{convolution} Let $\mo_1$ and $\mo_2$ be two
support-bounded constructible subsets of $ \mp^b(A,\bd'')$ and
$\mp^b(A,\bd'),$ respectively. Then we define
$$
\hat{1}_{\mo_{1}}*\hat{1}_{\mo_{2}}(L)=\chi(V(\mo_1,\mo_2;L))
$$
for $L\in\mp^b(A,\bd'+\bd'')$ and set
$F_{\mo_1,\mo_2}^{L}=\chi(V(\mo_1,\mo_2;L)).$
\end{Definition}

As Corollary \ref{constructible-well-defined}, we have
\begin{Prop}
If $\hat{f}\in M(\mp^b(A,\bd''))$ and $\hat{g}\in M(\mp^b(A,\bd')),$
then $\hat{f}*\hat{g}\in M(\mp^b(A,\bd'+\bd''))$ is well-defined.
\end{Prop}

Let $\hat{f}\in M(\mp^b(A,\bd''))$ and $\hat{g}\in
M(\mp^b(A,\bd')).$  Then
$$(\hat{f}*\hat{g})(x)=\int_{\ol{U}}\hat{f}(x')\hat{g}(x''):=\sum_{c,d\in\bbc}\chi(V(f^{-1}(c),g^{-1}(d);x)cd.$$
where $\ol{U}=V(\supp(f),\supp(g);x).$

\begin{Remark}
The definition \ref{convolution} does not supply an associative
multiplication in general. In \cite{Toen2005}, the author define an
associative multiplication for the derived category $\md^b(A)$ over
a finite field. However, it is not known how to make an analogy of
this associative multiplication over $\bbc.$
\end{Remark}

\section{The Relative homotopy category of $m$-cycle complexes}
Let $A$ be finite dimensional and finite global dimensional
associative algebra over $\mathbb{C}$ and $m$ be a positive even
number. We will recall some results in \cite{PengXiao1997} for the
relative homotopy category of $m$-cycle complexes over $A$ and
define their geometrization.
\subsection{}
A $m$-cycle complex over $A$ is by definition a complex
$X^{\blt}=(X_i,\partial_i)$ satisfying $X_{i}=X_j$ and
$\partial_i=\partial_j$ for all $i,j\in \mathbb{Z}$ with $i\equiv
j(mod\ m)$. If $X^{\blt}$ and $Y^{\blt}$ are two $m$-cycle
complexes, a \emph{morphism} $f^{\blt}:X^{\blt}\rightarrow
Y^{\blt}$ is a morphism of complexes such that $f_{i}=f_{j}$ for
all $i,j\in \mathbb{Z}$ with $i\equiv j(mod\ m)$. Hence, all
$m$-cycle complexes constitute an abelian subcategory of
$\mathcal{C}(A)$, denoted by $\mathcal{C}_{m}(A)$.We also denote
the subcategory of $\mathcal{C}_{m}(A)$ consisting of $m$-cycle
complex whose term is projective $A$-module by $\mp_{m}(A).$

Let $f^{\blt},g^{\blt}:X^{\blt}\rightarrow Y^{\blt}$ be two
morphisms of $m$-cycle complexes. A \emph{relative homotopy}
$s^{\blt}$ from $f^{\blt}$ to $g^{\blt}$ is a homotopy map of
complex such that $s_{i}=s_{j}$ for all $i,j\in \mathbb{Z}$ with
$i\equiv j(mod\ m)$. Under this condition, $f^{\blt}$ and
$g^{\blt}$ are said to be relatively homotopic. Hence, we can form
an additive and homotopy category $\mathcal{K}_{m}(A).$ We use
this notation though it usually is not a subcategory of
$\mathcal{K}(A)$. We also denote $\mp(A)\cap \mc_{m}(A)$ and
$\mathcal{K}(\mp(A))\cap \mathcal{K}_{m}(A)$ by $\mp_{m}(A)$ and
$\mathcal{K}_{m}(\mp(A))$.

We define a functor $CF:\mc^b(A)\rightarrow \mc_{m}(A)$ as
follows. For $X^{\blt}=(X_i,\partial^{X}_i)\in \mc^b(A)$, set
$FX^{\blt}=((FX^{\blt})_{i},\partial^{FX}_i)$ where
$(FX^{\blt})_{i}=\oplus X_{i+tm}$ and
$\partial^{FX}_i=(\partial^{X}_{i,s,t})_{s,t\in \mathbb{Z}}$ such
that $\partial^{X}_{i,s,t}:X_{i+sm}\rightarrow X_{i+1+tm}$ with
$\partial^{X}_{i,s,t}=0$ for $s \neq t$ and
$\partial^{X}_{i,s,t}=\partial^{X}_i$. $CF$ induces a functor
$F:\mathcal{K}^b(\mp(A))\rightarrow \mathcal{K}_{m}(\mp(A))$.
\begin{theorem}
$\mathcal{K}_{m}(\mp(A))$ is a triangulated category with the
shift functor defined for complex category and the functor
$F:\mathcal{K}^{b}(\mp(A))\rightarrow \mathcal{K}_{m}(\mp(A))$ is
exact.
\end{theorem}

For the natural functors:
$$\mk^b(\mp(A))\ra \mk^b(A)\ra\md^b(A)$$
the functors $F$ induce the following commutative diagram

$$
\xymatrix{\mk^b(\mp(A))\ar[d]^{F} \ar[r] & \mk^b(A) \ar[d]^{F}
\ar[r]& \md^b(A)\ar[d]^{F}\\\mk_m(\mp(A)) \ar[r] & \mk_m(A) \ar[r]
&\md_m(A) }
$$
where $\md_m(A)$ is the $m$-periodic derived category of $\mod A.$ It
is known that $\mk_m(\mp(A))$ is a triangulated full subcategory
of $\mk_m(A)).$ Therefore we can regard $\mk_m(\mp(A))$ as a
triangulated full subcategory of $\md_m(A).$ We consider the
triangulated full subcategory $\mr_m(A)$ of $\md_m(A)$ (also of
$\mk_m(\mp(A))$) generated by the full subcategory
$F(\mk^b(\mp(A))$ (also by $\mk^b(A)$). In general, the functor $F$
is not dense. If $A$ is hereditary, then $\mr_m(A)=\mk_m(A).$ When
$m=2$ we call $\mr_2(A)$ the root category of $\mod A.$ By the
way, the Galois group associated with $F$ is the cyclic group
generated by $T^m.$

\subsection{}
A complex $C^{\blt}$ of $A$-modules is called period-2 (or
2-periodic) complex if it satisfies $C^{\blt}[2]=C^{\blt}.$ We can
simply write the sequence of dimension vector of $C^{\blt}$ as
$\udd(C^{\blt})=(\ud(C^0),\ud(C^1)).$  For a sequence of dimension
vectors $\udd=(\ud_{0},\ud_1)$ , we define $\mc_{2}(A,\udd)$ to be
the subset of
$$\bbe_{\ud_0}(Q,R)\times\bbe_{\ud_1}(Q,R)\times\hom_{\bbc}(V^{\ud_0},V^{\ud_{1}})\times\hom_{\bbc}(V^{\ud_1},V^{\ud_{0}})$$
which consists of elements $x=(x^0,x^1,\pt_0,\pt_1),$ where
$x^i\in\bbe_{\ud_i}(Q,R)$  and $M(x^i)$ are the corresponding
$A$-modules
 on the space $V^{\ud_i}$ for $i=0,1$ respectively; and $\pt_0\in\hom_{\bbc}(V^{\ud_0},V^{\ud_{1}})$
 is a $A$-module homomorphism from $M(x^0)$
to $M(x^1)$,$\pt_1\in\hom_{\bbc}(V^{\ud_1},V^{\ud_{0}})$ is a
$A$-module homomorphism from $M(x^1)$ to $M(x^0)$ with the property
$\pt_1 \pt_0 =0$ and $\pt_0 \pt_1 =0$. As in \cite{LinPeng2002}, a
2-periodic complex $C^{\blt}$ can be written as
$$
\xymatrix{C^0\ar@<1ex>[r]^{\partial_0}&C^1\ar@<1ex>[l]^{\partial_1}}
$$
with $\partial_0\partial_1=\partial_1\partial_0=0.$

The group $G_{\udd}(Q)= GL_{\ud_0}(Q)\times GL_{\ud_1}(Q)$ acts on
$\mc_{2}(A,\udd) $ via the conjugation action
$$(g_i)_i (x_i,\pt_i)_i=((x_i)^{g_i},g_{i-1}\pt_i g_i^{-1})_i$$
where the action $(x_i)^{g_i}$ was defined as in Section 1.1. For
$x\in \mc_{2}(A,\udd)$, we denote its corresponding complex by
$M^{\blt}(x)$. Therefore the orbits under the action correspond
bijectively to the isomorphism classes of complexes of period-2
$A$-modules.

Giving the complete set $P_1,P_2,\cdots,P_l$ of indecomposable
projective $A$-modules as in Section 1. Let
$P^{\blt}=(P^0,P^1,\pt_0,\pt_1)$ be a period-$2$ complex where
$P^0\cong \oplus_{j=1}^l e^0_jP_j$ and $P^1\cong \oplus_{j=1}^l
e^1_jP_j$. Let $\ue(P^i)=(e^i_j)$ for $i=0,1.$ Then
$\uee=(\ue(P^0),\ue(P^1))$ is the projective dimension sequence of
$P^{\blt}.$  Write $\udd(\uee)=(\udim P^0,\udim P^1)$. We define
$\mp_{2}(A,\uee)$ to be the subset of $\mc_2(A, \udd(\uee))$
consisting of elements $x=(x^0,x^1,\pt_0,\pt_1),$ where
$(\bbc^{\ud_i}, x^i)\cong P^i$ for $i=0,1.$

The algebraic group $G_{\udd(\uee)}(Q,R)$ acts on $\mp_{2}(A,\uee)$
by
$$(g_i)_i(\pt_i)_i=(g_{i-1}\pt_i g_i^{-1})_i.$$
For $x\in \mp_{2}(A,\uee)$, we denote the corresponding complex in
$\mp_{2}(A)$ by $P^{\blt}(x)$.

\subsection{}
 A complex
$X^{\blt}=(X_i,\pt_i)_i\in \mp_{2}(A)$ is called contractible if the
induced homological groups $H_i(X^{\blt})=\ker\pt_{i+1}/ \Im\pt_i=0$
for all $i=0,1.$ It is easy to see that any contractible projective
complex is isomorphic to a direct sum of shifted copies of complexes
of the form
\begin{equation}\label{e11}
  \begin{CD}
  \dots  @>f >>P @>0>> P @>f>> P @>0 >>P @>f >>\dots
 \end{CD}
 \end{equation}
with $P$ a projective $A$-module and $f$ an automorphism. We call an
element $x\in\mp_2(A,\uee)$ contractible if the corresponding
projective complex is contractible. As in Section 2.3, we can prove
the $2$-periodic versions for Lemma \ref{lem1}, Lemma \ref{lem2} and
Lemma \ref{lem3}.  Any complex in $\mp_{2}(A)$ is isomorphic to the
direct sum of minimal period-2 projective complex and a contractible
projective complex. Hence, we can define a direct system
$\{(\mp_2(A, \uee), t_{\uee\uee'})\mid \uee,\uee'\in
\udim^{-1}(\bd), \uee\leq \uee' \}$ for any $\bd\in K_0$ where $K_0$
is the Grothendieck group of the category $\mk_2(\mp(A))$. We note
that  all triangulated categories $\mk_2(A),$ $\mk_2(\mp(A)),$
$\md_2(A),$ $\mk^b(\mp(A))$ and $\md^b(A)$ have the same
Grothendieck groups. Then we define
$$\mp_{2}(A,\bd):=\varinjlim_{\uee\in \udim^{-1}(\bd)}\mp_2(A, \uee) \quad\mbox{ and }\quad G_{\bd}=\varinjlim_{\uee\in \udim^{-1}(\bd)}G_{\udd(\uee)}.$$
We also have the quotient space
$$\mq\mp_{2}(A,\bd)=\mp_{2}(A,\bd)/\sim,$$
where $x\sim y$ in $\mp_{2}(A,\bd)$ if and only if the corresponding
projective complexes $P^{\blt}(x)$ and $P^{\blt}(y)$ are
quasi-isomorphic, i.e., they are isomorphic in $\mk_2(\mp(A)).$  We
denote by $\mq\mp_2(A, \uee)$ and $\mq_{1}\mp_2(A,\bd)$ the orbit
spaces of $\mp_2(A, \uee)$ and $\mp_2(A, \bd)$ under the actions of
$G_{\udd(\uee)}$ and $G_{\bd},$ respectively. As Proposition
\ref{quotientspace}, we have the follow analogy.
\begin{Prop}With the above notations, we have
$$
\mq\mp_2(A,\bd)=\mq_{1}\mp_2(A,\bd)=\varinjlim_{\udd(\uee)\in\underline{\mathrm{dim}}^{-1}(\bd)}\mq\mp_2(A,\uee).
$$
\end{Prop}

\subsection{} We can define the 2-periodic versions of all notations in Section 3.4 by
substituting $\mp_2(A, \uee)$ and $\mp_2(A, \bd)$ for $\mp^b(A,
\uee)$ and $\mp^b(A, \bd)$, respectively. For $\bd\in K_0$ and
$\uee\in \udim^{-1}(\bd),$ let $\mo\subseteq\mp_2(A, \uee)$ be a
$G_{\uee}$-invariant constructible subset. Then $G_{\bd}.\mo$ be a
support-bounded constructible subset. We denote by $\hat{1}_{\mo}$
the $\bbc$-value function over $\mp^b(A, \bd)$ given by taking
values $1$ on  each point in $\mo$ and $0$ otherwise. A function
$\hat{f}$ on $\mp_2(A,\bd)$ is called a $G_{\bd}$-invariant
constructible function
 if $\hat{f}$ can be written as a sum of finite terms $\sum_{i}m_i \hat{1}_{\mo_i}$
 where $m_i\in \bbc$ and  $\mo_i$ is a support-bounded constructible subset of $\mp_2(A,\bd)$. We denote by $M(\mp_2(A,
\bd))$ the set of $G_{\bd}$-invariant constructible functions.
\begin{Definition}\label{2-convolution} Let $\mo_1$ and $\mo_2$ be two
support-bounded constructible subsets of $ \mp_2(A,\bd'')$ and
$\mp_2(A,\bd'),$ respectively. Then we define
$$
\hat{1}_{\mo_{1}}*\hat{1}_{\mo_{2}}(L)=\chi(V(\mo_1,\mo_2;L))
$$
 and set
$F_{\mo_1,\mo_2}^{L}=\chi(V(\mo_1,\mo_2;L))$ for
$L\in\mp_2(A,\bd'+\bd'')$ where $V(\mo_1, \mo_2; L)$ is a 2-periodic
analogue of $V(\mo_1, \mo_2; L)$ defined in Section 3.4.
\end{Definition}

As Corollary \ref{constructible-well-defined}, we have
\begin{Prop}
If $\hat{f}\in M(\mp_2(A,\bd''))$ and $\hat{g}\in M(\mp_2(A,\bd')),$
then $\hat{f}*\hat{g}\in M(\mp_2(A,\bd'+\bd''))$ is well-defined.
\end{Prop}

Let $\hat{f}\in M(\mp_2(A,\bd''))$ and $\hat{g}\in
M(\mp_2(A,\bd')).$ Then
$$(\hat{f}*\hat{g})(x)=\int_{\ol{U}}\hat{f}(x')\hat{g}(x''):=\sum_{c,d\in\bbc}\chi(V(f^{-1}(c),g^{-1}(d);x)cd.$$
where $\ol{U}=V(\supp(f),\supp(g);x).$

\section{Realization of Lie algebras}
Let $\bd_1,\bd_2$ and $\bd$ be in $K_0$ and
$\mo_{1}\in\mp_2(A,\bd_1)$, $\mo_2\in\mp_2(A,\bd_2)$ and
$\mo\in\mp_2(A,\bd)$ be support-bounded constructible sets. Then
$\mo_1=G_{\bd_1}.\mo_{\uee''}, \mo_2=G_{\bd_2}.\mo_{\uee'}$ and
$\mo=G_{\bd}.\mo_{\uee}$ for projective dimension sequences $\uee'',
\uee'$ and $\uee.$

A constructible set $\mo$ is called indecomposable if all points
in $\mo$ correspond to indecomposable objects in $\mk_2(\mp(A)).$

\subsection{} Let $\ml=\ml(\mo_1,\mo_2)$ $$=\{L\in\mk_2(\mp(A))\mid \
\mbox{there exists triangle}\ Y\ra L\ra X\ra Y[1]\ \mbox{with}\
X\in\mo_1, Y\in \mo_2\}$$ Then $\ml$ is a support-bounded
constructible set in $\mp_2(A, \bd)$ and
$\ml=G_{\bd}.\ml_{\uee'+\uee''}$ for some constructible subset
$\ml_{\uee'+\uee''}$ in $\mp_2(A, \uee'+\uee'')$. We will consider
the following quotient stacks.

(a) Let
$W(\mo_{\uee''},\mo_{\uee'},\ml_{\uee'+\uee''}):=\bigcup_{L\in
\ml_{\uee'+\uee''}}W(\mo_{\uee''},\mo_{\uee'},L)$. The action of
$G_{\uee''}\times G_{\uee'}$ on
$W(\mo_{\uee''},\mo_{\uee'},\ml_{\uee'+\uee''})$ is defined as in
Section 3.4. For $(a,c)\in G_{\uee''}\times G_{\uee'},$
$(L,(f,g,h))\in W(\mo_{\uee''},\mo_{\uee'};\ml_{\uee'+\uee''})$,
define
$$(a,c)\circ(L,(f,g,h))=(L,(fc^{-1},ag,(c[1])ha^{-1})).$$
The quotient stack is $V(\mo_{\uee''}, \mo_{\uee'};
\ml_{\uee''+\uee'}):=[W(\mo_{\uee''}, \mo_{\uee'};
\ml_{\uee''+\uee'})/G_{\uee''}\times G_{\uee'}].$ We denote by
$(L,(f,g,h)^{\wedge})$ the geometric point in $V(\mo_{\uee''},
\mo_{\uee'}; \ml_{\uee''+\uee'}) $ corresponding to $(L, (f,g,h)).$
Moreover, up to $1$-isomorphism, the quotient stack is independent
of the choice $\uee''$ and $\uee'$ as in Section 3.4. Hence, we also
denote it by
$V(\mo_1,\mo_2;\ml):=V(\mo_{1\uee''},\mo_{2\uee'};\ml_{\uee''+\uee'})$
and set $F_{\mo_1\mo_2}^{\ml}=\chi(V(\mo_1,\mo_2;\ml)).$

(b) Let $M\in \mp_2(\uee''+\uee'+\uee)$ and
$$W_{(\mo_{\uee''}\mo_{\uee'})\mo_{\uee}}^{\mathcal{L}_{\uee''+\uee'} M}:=\bigcup_{L\in
\mathcal{L}_{\uee''+\uee'}}W(\mo_{\uee''},\mo_{\uee'};L)\times
W(L,\mo_{\uee};M).
$$ Consider the action of $G_{\uee''}\times G_{\uee'+\uee''}\times
G_{\uee'}\times G_{\uee}$. For $(a,b,c,d)\in G_{\uee''}\times
G_{\uee'+\uee''}\times G_{\uee'}\times G_{\uee}$ and $(L,
(f,g,h),(l,m,n))\in
W_{(\mo_{\uee''}\mo_{\uee'})\mo_{\uee}}^{\mathcal{L}_{\uee''+\uee'}
M},$ define
$$(a,b,c,d)\circ(L, (f,g,h),(l,m,n))=(L',(bfc^{-1},agb^{-1},c[1]ha^{-1}),(ld^{-1},bm,d[1]nb^{-1})).$$
The quotient stack is denoted by
$(W_{(\mo_{\uee''}\mo_{\uee'})\mo_{\uee}}^{\mathcal{L}_{\uee''+\uee'}
M})^{\wedge}$. We denote by
$$(L,(f,g,h),(l,m,n))^{\wedge}=\{(b(L),(bfc^{-1},agb^{-1},c[1]ha^{-1}),(ld^{-1},bm,d[1]nb^{-1}))\mid
$$
$$(a,b,c,d)\in G_{\uee''}\times G_{\uee'+\uee''}\times
G_{\uee'}\times G_{\uee}\}$$ the geometric point corresponding to
$(L, (f,g,h), (l, m,n)).$ The quotient stack is also independent of
the choices of $\uee'', \uee'$ and $\uee.$ Depending on the
discussion in Section 3.4, we can denote it by
$(W_{(\mo_1\mo_2)\mo}^{\mathcal{L} M})^{\wedge}$ and set
$\chi_{(\mo_1\mo_2)\mo}^{\mathcal{L}
M}=\chi((W_{(\mo_1\mo_2)\mo}^{\mathcal{L} M})^{\wedge}).$

Dually, let
$$W_{\mo_{\uee''}(\mo_{\uee'}\mo_{\uee})}^{ M\mathcal{L'}_{\uee+\uee'}}:=\bigcup_{L'\in
\mathcal{L'}_{\uee+\uee'}}W(\mo_{\uee},\mo_{\uee'}; L')\times W(
\mo_{\uee''}, L'; M).
$$ There is also an action of $G_{\uee''}\times G_{\uee'+\uee''}\times
G_{\uee'}\times G_{\uee}$ as above. For $(a,b',c,d)\in
G_{\uee''}\times G_{\uee'+\uee''}\times G_{\uee'}\times G_{\uee}$
and $(L', (f',g',h'),(l',m',n'))\in
W_{\mo_{\uee''}(\mo_{\uee'}\mo_{\uee})}^{
M\mathcal{L'}_{\uee+\uee'}},$ define
$$(a,b',c,d)\circ(L', (f',g',h'),(l,m,n))$$$$=(L'',(b'l'd^{-1},
cm'b'^{-1},(d[1])n'c^{-1}),(f'b'^{-1},ag',(b'[1])h'a^{-1})).$$ The
quotient stack is denoted by
$(W_{\mo_{\uee''}(\mo_{\uee'}\mo_{\uee})}^{
M\mathcal{L'}_{\uee+\uee'}})^{\wedge}$. We denote by
$$(L',(f',g',h'),(l',m',n'))^{\wedge}$$ the geometric point
corresponding to $(L', (f',g',h'), (l', m',n')).$ The quotient stack
is also independent of the choices of $\uee'', \uee'$ and $\uee.$
Depending on the discussion in Section 3.4, we can denote it by
$(W_{\mo_{\uee''}(\mo_{\uee'}\mo_{\uee})}^{
M\mathcal{L'}_{\uee+\uee'}})^{\wedge}$ and set
$\chi_{\mo_1(\mo_2\mo)}^{M\mathcal{L'}
}=\chi((W_{\mo_{\uee''}(\mo_{\uee'}\mo_{\uee})}^{
M\mathcal{L'}_{\uee+\uee'}})^{\wedge}).$

(c) The group $G_{\uee}$ acts on
$W(\ml_{\uee''+\uee'},\mo_{\uee};M)$ as follows. For any $d\in
G_{\uee},$ $(l,m,n)\in W(\ml_{\uee''+\uee'},\mo_{\uee};M),$
$$d\circ (l,m,n)=(ld^{-1},m,d[1]n)$$
The quotient stack is denoted by
$W(\ml_{\uee''+\uee'},\mo_{\uee};M)^{*}$.  It is independent of the
choice of $\uee$. Then we can also write $W(\ml_{\uee''+\uee'},\mo;
M)^{*}.$ We denote by $(l,m,n)^{*}=\{ld^{-1},m,d[1]n)\mid d\in
G_{\uee}\}$ the geometric point corresponding to $(l,m,n).$

\subsection {} Let $K_0$ be the Grothendieck group
of $\mk_2(\mp(A)).$ Write $\fh\simeq K_0\otimes_{\bbz}\bbc,$ which
is spanned by $\{h_{\bd}\mid\bd\in K_0\}$ subject to the relation:
$h_{\bd}=h_{\bd_1}+h_{\bd_2}$ if $\bd=\bd_1+\bd_2$ in $K_0.$ For
$\mo\subset\mp_2(A,\bd),$ we write $h_{\mo}:=h_{\bd}.$ The symmetric
Euler bilinear form on $\fh$ is given as
$$(h_{\bd_1}|h_{\bd_2})=\dim_{\bbc}\hom(X,Y)-\dim_{\bbc}\hom(X,Y[1])$$ $$ +\dim_{\bbc}\hom(Y,X)-\dim_{\bbc}\hom(Y,X[1])$$
for any $X\in\mp_2(A,\bd_1),Y\in\mp_2(A,\bd_2).$ It is well-defined.

Let $M(\mp_2(A,\bd))$ is the space of $G_{\bd}$-invariant
constructible functions over $\mp_2(A, \bd)$. A function $\hat{f}\in
M(\mp_2(A,\bd))$ is called indecomposable if any point in $supp(f)$
is indecomposable in $\mk_2(\mp(A))$. Let $I_{GT}(\bd)$ be the
$\bbc$-space of all indecomposable $G_{\bd}$-invariant functions in
$M(\mp_2(A,\bd)).$ We set
$$\fn=\bigoplus_{\bd\in K_0}I_{GT}(\bd)$$ and
$$\fg=\fh\oplus\fn.$$
We define the Lie bracket operation on $\fg=\fh\oplus\fn$ by the
following formulae:
$$
[\hat{1}_{\mo_{1}},\hat{1}_{\mo_{2}}]
=[\hat{1}_{\mo_{1}},\hat{1}_{\mo_{2}}]_{\fn}+\chi(\ol{\mo_1\cap\mo_2[1]})h_{\bd_1}
$$ where $\ol{\mo_1\cap\mo_2[1]}=
(\mo_{\uee''}\cap\mo_{\uee'}[1])/G_{\uee''}$ and $
[\hat{1}_{\mo_{1}},\hat{1}_{\mo_{2}}]_{\mathfrak n}(L):
=F^{L}_{\mo_1\mo_2}-F^{L}_{\mo_2\mo_1} $,
$$
[h_{\bd_1},\hat{1}_{\mo_2}]:= (h_{\bd_1}\mid
h_{\bd_2})\hat{1}_{\mo_2},
[\hat{1}_{\mo_2},h_{\bd_1}]=-(h_{\bd_1}\mid
h_{\bd_2})\hat{1}_{\mo_2}, $$
$$
[h_{\bd_1},h_{\bd_2}]:= 0.
$$
For $\hat{f}\in I_{GT}(\bd_1)$ and $\hat{g}\in I_{GT}(\bd_2),$ we
can write the formulae in the following integral form:
$$[\hat{f},\hat{g}]_{\fn}(x)=\int_{V(\supp(f),\supp(g);x)}f(x')g(x'')-\int_{V(\supp(g),\supp(f);x)}g(x')f(x'')$$
$$[\hat{f},h_{\bd}]=\int_{\supp(f)}f(x)(h_{\bd_1}|h_{\bd})\hat{1}_{f^{-1}(f(x))}.$$

We will prove the following two of the main theorems in this
paper.

\bigskip

\nd{\bf Theorem C} {\em Under the Lie bracket $[-,-]$ defined as
above, $\fg=\fh\oplus\fn$ is a Lie algebra over $\bbc.$}

\bigskip

\nd{\bf Theorem D} {\em There exists a symmetric bilinear form
$(-|-)$ defined on $\fg$ satisfying that
\begin{enumerate}
\item[(a)] $(-\mid-)$ is invariant in the following sense:
\begin{enumerate}
\item[(i)]$([h_{\bd},\hat{1}_{\mo_1}]\mid
\hat{1}_{\mo_2})=(h_{\bd}\mid [\hat{1}_{\mo_1},\hat{1}_{\mo_2}]);$

\item[(ii)]$([\hat{1}_{\mo_1},\hat{1}_{\mo_2}]\mid
\hat{1}_{\mo_3})=(\hat{1}_{\mo_1}\mid
[\hat{1}_{\mo_2},\hat{1}_{\mo_3}]).$
\end{enumerate}

\item[(b)]$(-\mid -)\mid_{\mathfrak h}$ is defined as above.
\item[(c)]$(\mathfrak h\mid \mathfrak n)=0.$

\item[(d)]for any $\mo_1,\mo_{2}$ indecomposable
$$
(\hat{1}_{\mo_1}\mid \hat{1}_{\mo_2})=\chi(\ol{\mo_{1}\cap
\mo_{2}[1]}).
$$
\item[(e)]$(-|-)|_{\mathfrak n\times\mathfrak n}$ is
non-degenerated.
\end{enumerate}}

\subsection{}
In this subsection, we compare the na\"\i ve Euler characteristics
of quotient stacks induced by that defined in Section 5.1. We use
the notations in Section 5.1. For fixed $L\in \ml_{\uee''+\uee'},$
we set $$\lr{L}:=\{E\in \ml_{\uee''+\uee'}\mid
F_{\mo_1\mo_2}^E=F_{\mo_1\mo_2}^L, F_{\mo_2\mo_1}^E=F_{\mo_2\mo_1}^L
\mbox{ and } \chi_{(\mo_1\mo_2)\mo}^{LM}=\chi_{(\mo_1\mo_2)\mo}^{E
M}\}.$$ In the same way as the proof of Proposition
\ref{constructible}, we know $\lr{L}$ is a constructible subset of
$\ml_{\uee''+\uee'}$ and there exists a finite subset
$R(\uee''+\uee')$ of $\ml_{\uee''+\uee'}$ such that
$$
\ml_{\uee''+\uee'}=\bigcup_{L\in R(\uee''+\uee')}\lr{L}.
$$
Then by Lemma \ref{Eulersum} and \ref{Eulermulti}, we obtain the
following result.

\begin{Lemma}
With the above notations, we have
$$
\chi^{\ml M}_{(\mo_1\mo_2)\mo}=\sum_{L\in
R(\uee''+\uee')}\chi^{LM}_{(\mo_1\mo_2)\mo}\cdot
\chi(\lr{L}/G_{\uee''+\uee'}).
$$
\end{Lemma}
Dually, for fixed $L'\in \ml_{\uee+\uee'},$ we set
$$\lr{L'}^{*}:=\{E\in \ml_{\uee+\uee'}\mid
F_{\mo_2\mo}^E=F_{\mo_2\mo}^{L'}, F_{\mo\mo_2}^E=F_{\mo\mo_2}^{L'}
\mbox{ and } \chi_{\mo_1(\mo_2\mo)}^{ML'}=\chi_{\mo_1(\mo_2\mo)}^{
ME}\}.$$ Then there exists a finite subset $R^*(\uee+\uee')$ of
$\ml'_{\uee+\uee'}$ such that
$$
\ml'_{\uee+\uee'}=\bigcup_{L\in R(\uee+\uee')}\lr{L'}^{*}.
$$

\begin{Lemma}
With the above notations, we have
$$
\chi^{ M\ml'}_{\mo_1(\mo_2\mo)}=\sum_{L'\in
R^*(\uee+\uee')}\chi^{ML'}_{(\mo_1\mo_2)\mo}\cdot
\chi(\lr{L}^*/G_{\uee+\uee'}).
$$
\end{Lemma}

\begin{Prop}
If $\mo_1,\mo_2$ are indecomposable, then
$$
\chi^{\ml M}_{(\mo_1\mo_2)\mo}={\chi^{M\ml' }_{\mo_1(\mo_2\mo)}}.
$$
\end{Prop}
Our proof follows \cite{PengXiao2000} and some improvements in
\cite{Hubery2005}.
\bigskip
\begin{proof}
  First we construct a map:
$$
W(\mo_{1\uee''},\mo_{2\uee'};\ml_{\uee''+\uee'})\times
W(\ml_{\uee''+\uee'},\mo_{\uee};M)\xrightarrow{\tau}
W(\mo_{1\uee''},\ml'_{\uee+\uee'};M)\times
W(\mo_{2\uee'},\mo_{\uee};\ml'_{\uee+\uee'})
$$
sending $((f,g,h),(l,m,n))$ to $((f',g',h'),(l',m',n'))$ via the
following commutative diagram:

\begin{equation}\label{327}
\xymatrix{Z\ar@{=}[r]\ar[d]^{l'} &Z\ar[d]^l\\
L'\ar@{.>}[r]^{f'}\ar[d]^{m'} &M\ar@{.>}[r]^{g'}\ar[d]^m &X\ar@{.>}[r]^-{h'}\ar@{=}[d] &L'[1]\ar[d]^{m'[1]}\\
Y\ar[r]^f\ar[d]^{n'} &L\ar[r]^g\ar[d]^n &X\ar[r]^-{h} &Y[1]\\
Z[1]\ar@{=}[r] &Z[1]}
\end{equation}

\nd where $L'$ is the direct summand of $Cone((m,f)[-1])$ such
that the complement is a contractible complex and $L'\in
\ml'_{\uee+\uee'}$ and we denote the natural embedding
$i:L'\rightarrow Cone((m,f)[-1])$. Hence $f'=(0,1,0)^{T}\circ i$
and $ m'=(0,0,1)^{T}\circ i,$ where $(0,1,0)^{T}$ and
$(0,0,1)^{T}$ are the natural projections:
$$
Cone((m,f)[-1])\rightarrow M\oplus Y\rightarrow M \mbox{  and
}Cone((m,f)[-1])\rightarrow M\oplus Y\rightarrow Y.
$$
Using the property of cone in the triangulated category,
$g',h',l',n'$ can be expressed by the composition of $f,g,h,l,m,n$
and some quasi-isomorphisms, by Lemma 2.1, 2.2 and 2.3, we deduce
that these expressions are algebraic. Hence, $\tau$ is a morphism
of varieties.

Moreover, we claim that $\tau$ is invariant under group action.
For $(a,b,c,d)\in G_{\uee''}\times G_{\uee'}\times G_{\uee},$
considering the pair $((\overline f,\overline g,\overline
h),(\overline l,\overline m,\overline
n))=(a,b,c,d)\circ((f,g,h),(l,m,n))$.
$$
\xymatrix{Y\ar[d]^{c}\ar[r]^{f} &L\ar[d]^{b}\ar[r]^{g}
&X\ar[d]^{a}\ar[r]^{h} &Y[1]\ar[d]^{c[1]}\\Y^{*}\ar[r]^{\overline f}
&L^{*}\ar[r]^{\overline g} &X^{*}\ar[r]^{\overline h} &Y^{*}[1]}
\qquad \xymatrix{Z\ar[d]^{d}\ar[r]^{l} &M\ar@{=}[d]\ar[r]^{m}
&L\ar[d]^{b}\ar[r]^{n} &Z[1]\ar[d]^{d[1]}\\Z^{*}\ar[r]^{\overline l}
&M\ar[r]^{\overline m} &L^{*}\ar[r]^{\overline n} &Z^{*}[1]}$$

Using the similar construction as diagram (\ref{327}), we again find
some $L''\in\ml'_{\uee+\uee'}$ and a pair $((\overline f',\overline
g',\overline h'),(\overline l',\overline m',\overline n'))\in
W(X^{*},L'';M)\times W(Y^{*},Z^{*};L'')$ satisfying:
$$
\tau((\overline f,\overline g,\overline h),(\overline l,\overline
m,\overline n))=((\overline f',\overline g',\overline h'),(\overline
l',\overline m',\overline n'))
$$
Following the construction of $L',$  there exists $b'$ such that
the following diagram of exact triangles commutes:
$$\xymatrix@C=15pt{L'\ar[rr]^-{(f'\,-m')^t}\ar[d]^{b'} &&M\oplus Y\ar[rr]^-{(m\,f)}\ar[d]^-{\left(\begin{smallmatrix}1\\&c\end{smallmatrix}\right)} &&L\ar[rr]^{h'g}\ar[d]^{b} &&L'[1]\ar[d]^{b'[1]}\\
L''\ar[rr]^-{(\overline f'\,-\overline m')^t} &&M\oplus
Y^{*}\ar[rr]^-{(\overline m\,\overline f)}
&&L^{*}\ar[rr]^{\overline{h}'\overline{g}} &&L''[1]}$$ It follows
that $b'$ is an isomorphism, and $b'\in G_{\uee+\uee'}.$

Similarly there exist $a'$ and $c'$ giving the following commutative
diagram of exact triangles:
$$\xymatrix{L'\ar[d]^{b'}\ar[r]^{f'} &M\ar@{=}[d]\ar[r]^{g'} &X\ar[d]^{a'}\ar[r]^{h'} &L'[1]\ar[d]^{b'[1]}\\L''\ar[r]^{\overline f'} &M\ar[r]^{\overline g'} &X^{*}\ar[r]^{\overline h'} &L''[1]} \qquad
\xymatrix{Z\ar[d]^{d'}\ar[r]^{l'} &L'\ar[d]^{b'}\ar[r]^{m'}
&Y\ar[d]^{c}\ar[r]^{n'} &Z[1]\ar[d]^{d'[1]}\\Z^{*}\ar[r]^{\overline
l'} &L'\ar[r]^{\overline m'} &Y^{*}\ar[r]^{\overline n'}
&Z^{*}[1]}$$ This shows $((\overline f',\overline g',\overline
h'),(\overline l',\overline m',\overline
n'))=(a',b',c,d')\circ((f',g',h'),(l',m',n'))$ and so the two pairs
of exact triangles lie in the same orbit. Hence, The morphism $\tau$
induces the following 1-morphism:
$$
\xymatrix{ (W(\mo_{1},\mo_{2};\ml)\times
W(\ml,\mo;M))^{\wedge}\ar[r]^{{\tau}^{\wedge}}&
(W(\mo_{1},\ml';M)\times W(\mo_{2},\mo;\ml'))^{\wedge} }
$$
Then we have a pseudomorphism $$\xymatrix{
(W(\mo_{1},\mo_{2};\ml)\times
W(\ml,\mo;M))^{\wedge}\ar[r]^{({\tau}^{\wedge})_*}(\bbc)&
(W(\mo_{1},\ml';M)\times W(\mo_{2},\mo;\ml'))^{\wedge} }(\bbc).$$
Depending on the symmetry of the construction of $L'$ , we can
construct the inverse of $({\tau}^{\wedge})_*$. Therefore,
$({\tau}^{\wedge})_*$ is a pseudoisomorphism.  The proof is
completed.
\end{proof}

We introduce some notations. Let
$$W_{(\mo_{\uee''}\mo_{\uee'})\mo_{\uee}}^{\mathcal{L}_{\uee''+\uee'} M}:=\bigcup_{L\in
\mathcal{L}_{\uee''+\uee'}}W(\mo_{\uee''},\mo_{\uee'};L)\times
W(L,\mo_{\uee};M)
$$ and
let $c(L)=((f,g,h)(l,m,n))$ denotes the triangles as follows:
\begin{equation}\label{tri}
Y\xrightarrow{f}L\xrightarrow{g}X\xrightarrow{h}Y[1],
\hspace{1cm}Z\xrightarrow{l}M\xrightarrow{m}L\xrightarrow{n}Z[1].
\end{equation}
We define
$$
\hspace{-0.9cm}W_{(\mo_{\uee''}\mo_{\uee'})\mo_{\uee}}^{\ml_{\uee''+\uee'}M}(1)=\{c(L)\in
W_{(\mo_{\uee''}\mo_{\uee'})\mo_{\uee}}^{\ml_{\uee''+\uee'}M}\mid
L\ncong M\oplus Z[1] \mbox{ for any } Z\in\mo_{\uee} \},
$$
$$
W_{(\mo_{\uee''}\mo_{\uee'})\mo_{\uee}}^{\ml_{\uee''+\uee'}M}(2)=\{c(L)\in
W_{(\mo_{\uee''}\mo_{\uee'})\mo_{\uee}}^{\ml_{\uee''+\uee'}M}\mid
L\cong M\oplus Z[1],  L\ncong X\oplus Y\mbox{ and } X\ncong Y\},
$$
$$
W_{(\mo_{\uee''}\mo_{\uee'})\mo_{\uee}}^{\ml_{\uee''+\uee'}M}(3)=\{c(L)\in
W_{(\mo_{\uee''}\mo_{\uee'})\mo_{\uee}}^{\ml_{\uee''+\uee'}M}\mid
L\cong M\oplus Z[1],  L\cong X\oplus Y\mbox{ and } X\ncong Y\},
$$
and
$$
W_{(\mo_{\uee''}\mo_{\uee'})\mo_{\uee}}^{\ml_{\uee''+\uee'}M}(4)=\{c(L)\in
W_{(\mo_{\uee''}\mo_{\uee'})\mo_{\uee}}^{\ml_{\uee''+\uee'}M}\mid
L\cong M\oplus Z[1],  L\cong X\oplus Y\mbox{ and } X\cong Y\}.
$$
It is clear that
$$
W_{(\mo_{\uee''}\mo_{\uee'})\mo_{\uee}}^{\ml_{\uee''+\uee'}M}=\bigcup_{i=1}^{4}W_{(\mo_{\uee''}\mo_{\uee'})\mo_{\uee}}^{\ml_{\uee''+\uee'}M}(i).
$$
The above partition induces the partitions of
$(W_{(\mo_{1}\mo_{2})\mo}^{\ml M})^{\wedge}.$ Hence,
$$
\chi_{(\mo_1\mo_2)\mo}^{\ml
M}=\sum_{i=1}^{4}\chi_{(\mo_1\mo_2)\mo}^{\ml M}(i)
$$
where $\chi_{(\mo_1\mo_2)\mo}^{\ml
M}(i)=\chi((W_{(\mo_{1}\mo_{2})\mo}^{\ml M})^{\wedge}(i)).$

Define
$$
V^{\ml M}_{(\mo_1\mo_2)\mo}=\bigcup_{L\in R(\uee''+\uee')}V(\mo_1,
\mo_2; L)\times V(\lr{L}, \mo; M)
$$

 Put $\lr{L}_0=\{E\in \lr{L}\mid E\cong M\oplus Z[1] \mbox {
for some } Z\in \mo_{\uee}\}$ and $\lr{L}_1=\lr{L}\backslash
\lr{L}_0.$ Hence,
$$
V^{\ml M}_{(\mo_1\mo_2)\mo}=V^{\ml M}_{(\mo_1\mo_2)\mo}(0)\bigcup
V^{\ml M}_{(\mo_1\mo_2)\mo}(1)
$$
where $V^{\ml M}_{(\mo_1\mo_2)\mo}(i)=\bigcup_{L\in
R(\uee''+\uee')}V(\mo_1, \mo_2; L)\times V(\lr{L}_i, \mo; M)$ for
$i=0,1.$
\begin{Lemma}\label{eulerconstant}
For fixed $M$ and $L\in \ml_{\uee''+\uee'},$ we have
$$
\chi(V(\mo_1,\mo_2;L)\times V(\mo,\lr{L}_1;
M))=\chi_{(\mo_1\mo_2)\mo}^{\lr{L}_1 M}.
$$
\end{Lemma}
\begin{proof}
Consider the $G_{\uee''+\uee'}$-action $$W(\mo_{\uee}, \lr{L}_1;
M)^*\rightarrow V(\mo, \lr{L}_1; M)
$$
For any $(l,m,n)^{*}$, the stable subgroup denoted by $B(n)$ is
isomorphic to the affine space $\mathrm{Hom}(Z[1], M)n.$ Then we
have a $1$-morphism
$$
p: (W_{(\mo_{\uee''}\mo_{\uee'})\mo_{\uee}}^{\lr{L}_1
M})^{\wedge}\rightarrow V(\mo, \lr{L}_1; M).
$$
It induces the pseudomorphism $$p_*:
(W_{(\mo_{\uee''}\mo_{\uee'})\mo_{\uee}}^{\lr{L}_1
M})^{\wedge}(\bbc)\rightarrow V(\mo, \lr{L}_1; M)(\bbc).$$ For any
$(l,m,n)^{\wedge}\in V(\mo, \lr{L}_1; M),$ the fibre
$p_*^{-1}(l,m,n)$ is pseudomorphic to
$[V(\mo_{1\uee''},\mo_{2\uee'}; E)/B(n)](\bbc)$ where $E$ occurs in
the triangle $(l,m,n)$ as diagram \eqref{tri}. Under the action of
$B(n),$ the stable subgroup for any $(f,g,h)\in
V(\mo_{1\uee''},\mo_{2\uee'}; \lr{L}_1)$ is
$$
\{b\in B(n)\mid fa'=bf \mbox{ for some } a'\in \mathrm{End}Y\}
$$
which is the subspace of $\mathrm{Hom}(Z[1],L)n.$ By Proposition
\ref{Euler}, Lemma \ref{Eulermulti} and  \ref{Eulersum}, according
to the fact $F_{\mo_1\mo_2}^E=F_{\mo_1\mo_2}^L$ for any $E\in
\lr{L}_1,$ we have
$$ \chi(V(\mo_1,\mo_2;L)\times V(\mo,\lr{L}_1;
M))=\chi_{(\mo_1\mo_2)\mo}^{\lr{L}_1 M}.
$$
\end{proof}
The Lemma naturally induces the following Proposition.
\begin{Prop}\label{identity1}
For fixed $M,$ we have
$$\chi_{(\mo_1\mo_2)\mo}^{\ml M}(1)=V^{\ml M}_{(\mo_1\mo_2)\mo}(1).$$
\end{Prop}

The following Lemma is a natural corollary of Proposition
\ref{Euler}.
\begin{Lemma}
Let $X,Y\in \mp_2(A)$ be two indecomposable objects, then
$$\chi(\mathrm{Hom}_{\mp_2(A)}(X,Y))=1 \mbox{ and }\chi(\mathrm{Aut}_{\mp_2(A)}X)=0.$$
\end{Lemma}

\begin{Prop}\label{number}
Let $\mo_1,\mo_2,\mo$ be indecomposable as above.
\begin{enumerate}
\item[(I)] If $L\cong M\oplus Z[1]$ for some $Z\in \mo$ and
$L\notin \mo_1\oplus\mo_2$, then $F_{\mo_1\mo_2}^{L}=0$ and
 $\chi_{(\mo_1\mo_2)\mo}^{\ml
M}(2)=\chi_{\mo(\mo_1\mo_2)}^{M\ml}(2).$ \item[(II)]If $X\in
\mo_1,Y\in \mo_2$ such that $X\ncong Y$, then
$\chi_{(XY)Y[1]}^{X\oplus Y, X}=\chi_{X[1](XY)}^{Y
 ,X\oplus Y}=1,\ \chi_{(XY)X[1]}^{X\oplus Y,Y }-1=\chi_{Y[1](XY)}^{X,X\oplus Y}-1=\mathrm{dim}_{\mathbb{C}}\mathrm{Hom}(Y,X).$
\item[(III)] If $X$ is indecomposable, then
$\chi_{(XX)X[1]}^{X\oplus X,X}=\chi_{X[1](XX)}^{X,X\oplus X}.$
Hence, $\chi_{(\mo_1\mo_2)\mo}^{\ml
M}(4)=\chi_{\mo(\mo_1\mo_2)}^{M\ml}(4).$
\end{enumerate}
\end{Prop}
\begin{proof}   Let $L\in M\oplus \mo[1]$.  Then
$L\cong M\bigoplus Z[1]$ for some $Z\in \mo.$ Define
$$
\phi: V(\mo_1,\mo_2;M\oplus Z[1])\rightarrow
((W(\mo_1,\mo_2;M\oplus Z[1])\times W(M\oplus Z[1],Z;M))^{\wedge}
$$
mapping $(f,g,h)^{\wedge}$ to $((f,g,h),(0, (1,0)^t,
(0,1)))^{\wedge}.$ It is a surjective morphism.  The group action
$W(M\oplus Z[1],Z;M)^*\rightarrow V(M\oplus Z[1],Z;M)$ has the
stable subgroup $B(n)$. It is independent of $n$ since $V(M\oplus
Z[1], Z; M)$ has only one element. We denote it by $B.$ It is
given as follows.
$$
B=\{b=\left(
        \begin{array}{cc}
          1 & 0 \\
          b_1 & b_2 \\
        \end{array}
      \right)\in \mathrm{End}(M\oplus Z[1])\mid b_{1}\in \mathrm{Hom}(Z[1],M),b_{2}\in \mathrm{Aut}
Z[1] \}
$$
Then $Stab.(f,g,h)^{\wedge}=\{b\in B\mid fc=bf $ for some $c\in
\mathrm{Aut}Y\}.$

For the case $L\notin \mo_1\bigoplus \mo_2,$ then
$$
1-Stab.((f,g,h)^{\wedge})=\{b'=\left(
        \begin{array}{cc}
          0 & 0 \\
          b'_1 & b'_2 \\
        \end{array}
      \right)\mid b'_1\in \mathrm{Hom}(Z[1],M), b'_2\in \mathrm{radEnd}Z[1]$$$$\hspace{-1cm} \mbox{ such that }fc=b'f \mbox{ for some }c\in
\mathrm{End}Y\}
$$
(see Section 7.3 of \cite{PengXiao2000}). It is an affine space.

 For any $x \in ((W(\mo_1,\mo_2;M\oplus Z[1])\times
W(M\oplus Z[1],Z;M))^{\wedge},$  by  Lemma \ref{Eulermulti} and
Lemma \ref{eulerconstant},
$$
\chi(\phi^{-1}(x))=\chi(B)=\chi(\mathrm{Hom}(Z[1],M))\cdot
\chi(\mathrm{Aut}Z[1])=0
$$
\\
and clearly $F_{M\oplus Z[1],Z}^{M}=1.$  Therefore
$$
F_{\mo_1\mo_2}^{M\oplus Z[1]}=\chi(V(\mo_1,\mo_2;M\oplus Z[1]))=0.
$$
The maps $\phi$ and $\varphi$ induce a homeomorphism:
$$
\hspace{-5cm}((W(\mo_1,\mo_2;M\oplus Z[1])\times W(Z,M\oplus
Z[1];M))^{\wedge}$$
$$\hspace{3cm}\rightarrow ((W(\mo_1,\mo_2;M\oplus Z[1])\times W(M\oplus
Z[1],Z;M))^{\wedge}
$$
Hence, $\chi_{(\mo_1\mo_{2})Z}^{M\oplus
Z[1],M}=\chi_{Z(\mo_1\mo_{2})}^{M,M\oplus Z[1]}.$ The statement
(I) is proved.  The statement (III) can be proved by the similar
discussion.

For the case  $L\in \mo_1\bigoplus \mo_2,$ then $L\cong X\oplus Y$
for some $X\in\mo_1, Y\in\mo_2$ and $X\ncong Y.$  Every orbit in
$V(X,Y,X\oplus Y)$ is of the form as
$(\left(%
\begin{array}{c}
  \theta_1 \\
  1 \\
\end{array}%
\right),(1,\theta_2),0)^{\wedge}$ such that
$\theta_1,\theta_2:Y\rightarrow X,\theta_1+\theta_2=0$. So
$V(X,Y;X\oplus Y)\cong \mathrm{Hom}(Y,X).$  This induces the maps:
$$
\mathrm{Hom}(Y,X)\rightarrow (W(X,Y;X\oplus Y)\times W(X\oplus
Y,Y[1];X))^{\wedge}
$$
and
$$
\mathrm{Hom}(Y,X)\rightarrow (W(X,Y;X\oplus Y)\times W(X\oplus
Y,X[1];Y))^{\wedge}
$$
The set $(W(X,Y;X\oplus Y)\times W(X\oplus Y,Y[1];X))^{\wedge}$
has unique element, implying $\chi_{(XY)Y[1]}^{X\oplus Y,X}=1$.
 The natural action of group $\mathrm{Aut}X$ on
$\mathrm{Hom}(Y,X)$ is free except on the point $0.$ It induces a
homeomorphism: $$(W(X,Y;X\oplus Y)\times W(X\oplus
Y,X[1];Y))^{\wedge}\setminus\{0\}\cong
(\mathrm{Hom}(Y,X)\setminus\{0\})/\mathrm{Aut}X$$  So
$\chi_{(XY)X[1]}^{X\oplus
Y,Y}=1+\chi(\mathrm{Hom}(Y,X)\setminus\{0\}/\mathrm{Aut}X)$. As we
know( see \cite{Riedtmann1994}),  $AutX$ is the direct product of
$\mathbb{C}^{*}$ and the subgroup $1+\mathrm{radEnd}X$ which is
contractible. Hence,
$$
\chi(\mathrm{Hom} (Y,X)\setminus\{0\}/\mathrm{Aut}X)= \chi(
\mathrm{Hom}(Y,X)\setminus\{0\}/\mathbb{C}^{*})=\chi(\mathbb{C}\mathbb{P}^{\mathrm{dim}_{\mathbb{C}}
\mathrm{Hom}(Y,X)})$$$$=\mathrm{dim}_{\mathbb{C}}\mathrm{Hom}(Y,X).
$$ The statement (II) is proved.
\end{proof}

\begin{Lemma}
If $L$ is decomposable or 0, then
$F_{\mo_1\mo_2}^L-F_{\mo_2\mo_1}^L= 0$.
\end{Lemma}
\begin{proof}  If $L$ is 0, it is trivial. If $L$ is decomposable, we assume
$L=M\oplus Z[1]$ with $Z$ indecomposable and $M\not\cong 0$. If
$L\not\cong X\oplus Y$ for any $X\in\mo_1,Y\in\mo_2$, then by
Proposition 5.5, $F_{\mo_1\mo_2}^L=F_{\mo_1\mo_2}^L=0.$ Using this
fact, we have:
$$
F_{\mo_1\mo_2}^L-F_{\mo_2\mo_1}^L=(F_{MZ[1]}^{M\oplus
Z[1]}-F_{Z[1]M}^{M\oplus
Z[1]})(1_{\mo_1}(M)1_{\mo_2}(Z[1])-1_{\mo_1}(Z[1])1_{\mo_2}(M))
$$
So we can assume that both $M,Z[1]$ are indecomposable and
$M\ncong Z[1]$. By the proof of II of Proposition \ref{number},
$F_{MZ[1]}^{M\oplus Z[1]}=\chi(\mathrm{Hom}(Z[1],M))=1$,
$F_{Z[1]M}^{M\oplus Z[1]}=\chi(\mathrm{Hom}(M,Z[1]))=1.$
\end{proof}
\begin{Cor}\label{number2}
For fixed $M$, we have
$$
\chi(V^{\ml M}_{(\mo_1\mo_2)\mo}(0))-\chi(V^{\ml
M}_{(\mo_2\mo_1)\mo}(0))=0.
$$
and
$$
\chi(V^{\ml M}_{\mo(\mo_1\mo_2)}(0))-\chi(V^{\ml
M}_{\mo(\mo_2\mo_1)}(0))=0.
$$
\end{Cor}
We note that we can translate the discussion in this section for
$W^{M\ml'}_{\mo_1(\mo_2\mo)}$ and $V^{M\ml' }_{\mo_1(\mo_2\mo)}$.
the discussion is completely dual.

\bigskip
\subsection{}
Now, we come to prove Theorem {\bf C} and {\bf D}.

\bigskip
\nd{\bf Proof of Theorem C:}  It is suffices to verify the Jacobi
identity:
$$
\Delta=[[\hat{1}_{\mo_1},\hat{1}_{\mo_2}],\hat{1}_{\mo}]-[[\hat{1}_{\mo_1},\hat{1}_{\mo}],\hat{1}_{\mo_2}]-[[\hat{1}_{\mo},\hat{1}_{\mo_2]},\hat{1}_{\mo_1}]=
0
$$
$$
\Delta'=[[h_{\mo_1},\hat{1}_{\mo_2}],\hat{1}_{\mo}]-[[h_{\mo_1},\hat{1}_{\mo}],\hat{1}_{\mo_2}]-[[\hat{1}_{\mo},\hat{1}_{\mo_2]},h_{\mo_1}]=
0
$$
$$
\Delta''=[[h_{\mo_1},h_{\mo_2}],\hat{1}_{\mo}]-[[h_{\mo_1},\hat{1}_{\mo}],h_{\mo_2}]-[[\hat{1}_{\mo},h_{\mo_2]},h_{\mo_1}]=
0
$$

\nd We will follow tightly  \cite{PengXiao2000} and
\cite{Hubery2005}. First, $\Delta$ (also $\Delta'$ and $\Delta''$)
is a $G_{\bd_1+\bd_2+\bd}$-invariant constructible function by
discussion in Section 3.4. Moreover, according to the definition of
Lie bracket operation in Section 5.2 and the discussion in Section
5.3, we know
$$
(\hat{1}_{\mo_1}*\hat{1}_{\mo_2})*\hat{1}_{\mo}(M)=\sum_{L\in
R(\uee''+\uee')}F_{\mo_1\mo_2}^{L}\cdot F_{\lr{L}\mo_1}^{M}$$ for
$M\in \mp_2(A,\bd_1+\bd_2+\bd)$ and
$$[[\hat{1}_{\mo_1},\hat{1}_{\mo_2}],\hat{1}_{\mo}]
=[[\hat{1}_{\mo_1},\hat{1}_{\mo_2}]_{\mathfrak
n},\hat{1}_{\mo}]_{\mathfrak n}-\chi(\overline{\mo_1\cap
\mo_{2}[1]})(h_{\bd_1}\mid
h_{\bd})\hat{1}_{\mo}-(F_{\mo_1\mo_2}^{\mo[1]}-F_{\mo_2\mo_1}^{\mo[1]})h_{\bd}
$$
where
$$[[\hat{1}_{\mo_1},\hat{1}_{\mo_2}]_{\mathfrak n},\hat{1}_{\mo}]_{\mathfrak
n}(M)=\sum_{L\in
R(\uee''+\uee')}(F_{\mo_1\mo_2}^{L}-F_{\mo_2\mo_1}^L)\cdot
(F_{\lr{L}\mo_1}^{M}-F_{\mo_1\lr{L}}^M)
$$
and $$F_{\mo_1\mo_2}^{\mo[1]}=\sum_{L\in
R(\uee''+\uee')}\chi(\lr{L}\cap \mo[1])F_{\mo_1\mo_2}^L.
$$ Let us write
$$c_M:=\Delta_{\mo_1\mo_2\mo}^M+\Delta_{\mo_2\mo\mo_1}^M+\Delta_{\mo\mo_1\mo_2}^M-\Delta_{\mo_2\mo_1\mo}^M-\Delta_{\mo\mo_2\mo_1}^M-\Delta_{\mo_1\mo\mo_2}^M,$$
where
$$\Delta_{\mo_1\mo_2\mo}^{M}:=\chi(V_{(\mo_1\mo_2)\mo}^{\ml M})-\chi(V_{\mo_1(\mo_2\mo)}^{M\ml}).$$
Define
$$
b_M:=(\chi(\overline{\mo_1\cap \mo_{2}[1]})(h_{\mo_1}\mid
h_{\mo})\hat{1}_{\mo}-\chi(\overline{\mo_1\cap
\mo[1]})(h_{\mo_1}\mid h_{\mo_2})\hat{1}_{\mo_2}$$
$$-\chi(\overline{\mo\cap \mo_{2}[1]})(h_{\mo}\mid
h_{\mo_1})\hat{1}_{\mo_1})(M).
$$
Set
$\gamma_{\mo_1\mo_2}^{\ml}=F_{\mo_1\mo_2}^{\ml}-F_{\mo_2\mo_1}^{\ml},$
we define
$$
a_M:=\gamma_{\mo_1\mo_2}^{\mo[1]}h_{\mo}
-\gamma_{\mo_1\mo}^{\mo_2[1]}h_{\mo_2}-\gamma_{\mo\mo_2}^{\mo_1[1]}h_{\mo_1}$$
hence,
$$
\Delta(M)=c_{M}-b_{M}-a_{M}
$$

\nd (1) $c_{M}=c^{1}_{M}+c^{2}_{M}$ and $c^{1}_{M}=0$

By Proposition 5.3 and $\ref{number}$, we have
$$ \Delta_{\mo_1\mo_2\mo}^{M}=
\chi(V_{(\mo_{1}\mo_{2})\mo}^{\ml
M}(0)))-\chi(V_{\mo_{1}(\mo_{2}\mo)}^{\ml M}(0)))
+\sum_{i=2}^4(\chi_{\mo_1(\mo_2\mo)}^{M\ml'}(i)-\chi_{(\mo_1\mo_2)\mo}^{\ml
M}(i)).
$$ Set
$$
c_{M}^{1} =(\chi(V_{(\mo_{1}\mo_{2})\mo}^{\ml M}(0))
-\chi(V_{(\mo_{2}\mo_{1})\mo}^{\ml M}(0)))
-(\chi(V_{\mo_1(\mo_{2}\mo)}^{M\ml }(0))
-\chi(V_{\mo_1(\mo\mo_{2})}^{M\ml }(0)))
$$$$\hspace{1cm}+(\chi(V_{(\mo_{2}\mo)\mo_1}^{\ml M }(0))
-\chi(V_{(\mo\mo_{2})\mo_1}^{\ml M}(0)))
-(\chi(V_{\mo_2(\mo\mo_{1})}^{M\ml }(0))
-\chi(V_{(\mo_{2}(\mo_{1}\mo)}^{M\ml
}(0)))$$$$\hspace{1cm}+(\chi(V_{(\mo\mo_1)\mo_{2})}^{\ml M}(0))
-\chi(V_{(\mo_{1}\mo)\mo_2}^{M\ml }(0))))
-(\chi(V_{\mo(\mo_1\mo_{2})}^{M\ml }(0))
-\chi(V_{(\mo(\mo_{2}\mo_{1})}^{M\ml }(0)))$$
$$c_{M}^{2}= \sum_{i=2}^4\{-\chi_{(\mo_1\mo_2)\mo}^{\ml
M}(i)+\chi_{\mo_1(\mo_2\mo)}^{M\ml'}(i)-\chi_{(\mo_2\mo)\mo_1}^{\ml
M}(i)+\chi_{\mo_2(\mo\mo_1)}^{M\ml'}(i)-\chi_{(\mo\mo_1)\mo_2}^{\ml
M}(i)$$$$\hspace{0.5cm}+\chi_{\mo(\mo_1\mo_2)}^{M\ml'}(i)
+\chi_{(\mo_2\mo_1)\mo}^{\ml
M}(i)-\chi_{\mo_2(\mo_1\mo)}^{M\ml'}(i)+\chi_{(\mo\mo_2)\mo_1}^{\ml
M}(i)-\chi_{\mo(\mo_2\mo_1)}^{M\ml'}(i)$$$$\hspace{-6cm}+\chi_{(\mo_1\mo)\mo_2}^{\ml
M}(i)-\chi_{\mo_1(\mo\mo_2)}^{M\ml'}(i)\}.$$

\nd Then $c_{M}=c_{M}^{1}+c_{M}^{2}$.
\\
By Corollary $\ref{number2}$, $c_M^{1}=0$ and for $c_M^{2}$, we
first remark the following fact (see  Proposition \ref{number}):
$$\sum_{t=2}^4\chi_{(\mo_i\mo_j)\mo_k}^{\ml
M}(t)-\chi_{\mo_k(\mo_i\mo_j)}^{M\ml'}(t)=
\chi_{(\mo_i\mo_j)\mo_k}^{\ml
M}(3)-\chi_{\mo_k(\mo_i\mo_j)}^{M\ml'}(3)$$
$$=\int_{Z\in\overline{\mo_{i}[1]\cap\mo_k}}\mathrm{dim}_{\mathbb{C}}\mathrm{Hom}(M,Z[1])1_{\mo_j}(M)-\int_{Z\in\overline{\mo_{j}[1]\cap\mo_k}}\mathrm{dim}_{\mathbb{C}}\mathrm{Hom}(Z[1],M)1_{\mo_i}(M)$$
for any $Z\in \mo_k.$ We use this fact to substitute the terms in
$c_{M}^{2}$, then
$$
c_{M}^{2}=b_{M}.
$$

\nd (2) $a_M=0$
\\
We claim that for indecomposable constructible sets $\mo_i,\mo_j$
and $\mo_k,$
\begin{equation}\label{1649}
F_{\mo_i\mo_j}^{\mo_k[1]}=F_{\mo_j\mo_k}^{\mo_i[1]}.
\end{equation}
 Without loss of
generality, we assume that $\mo_s$ is an indecomposable
constructible subset of $\mp_2(A,\uee_s)$ for $s=i,j,k.$
 Let $X_i,$
$X_j$ and $X_k$ be any objects in them, respectively and the
corresponding orbits are denoted by $\mo_{X_l}$ for $l=i,j,k.$
Considering group action $(c)$ in Section 5.1, we find as the
orbit spaces of $V(X_i,X_j;X_k[1])$ and $V(X_j[1],X_k[1];X_i)$
$$
\overline{V(X_i,X_j;X_k[1])}\cong \overline{V(X_j[1],X_k[1];X_i)}
$$
and the corresponding fibres are:
$$
\{(bf,gb^{-1},h)^{\wedge}\mid b\in G_{\uee}\}\mbox{ and
}\{(ag,ha^{-1},f[1])^{\wedge}\mid a\in G_{\uee''}\}
$$
Depending on indecomposable property of $X_i,X_j$ and $X_k,$,
$\mathrm{Aut}X_i$ and $\mathrm{Aut}X_k$ can be decomposed as direct
products of $\mathbb{C}^{*}$ and contractible subgroups, hence, we
can regard $a,b\in \mathbb{C}^{*}$ in the above fibres so that two
fibres have the same Euler characteristic. This shows
$F_{\mo_{X_i}\mo_{X_j}}^{X_k[1]}=F_{\mo_{X_j}\mo_{X_k}}^{X_i[1]}.$
Hence, by Lemma \ref{Eulermulti}, we know
$$
F_{\mo_{X_i}\mo_{X_j}}^{\mo_{X_k}[1]}=F_{\mo_{X_j}\mo_{X_k}}^{\mo_{X_i}[1]}.
$$

Let $\lr{X_i,X_j,X_k}=\{(E_i,E_j,E_k)\in \mo_i\times\mo_j\times
\mo_k \mid
F_{\mo_{E_i}\mo_{E_j}}^{\mo_{E_k}[1]}=F_{\mo_{X_i}\mo_{X_j}}^{\mo_{X_k}[1]}\}.$
Then
$$
\lr{X_i,X_j,X_k}=\lr{X_j,X_k,X_i}.
$$
Consider the projection
$$
W(\mo_i,\mo_j;\mo_k)\rightarrow \mo_i\times\mo_j\times\mo_k.
$$
It induces the map:
$$
V(\mo_i,\mo_j;\mo_k)\rightarrow
\overline{\mo}_i\times\overline{\mo}_j\times\overline{\mo}_k
$$
where $\overline{\mo_i\times\mo_j\times\mo_k}$ is the quotient
space of $\mo_i\times\mo_j\times\mo_k$ under the action of
$G_{\uee_i}\times G_{\uee_j}\times G_{\uee_k}$. In the same way as
the proof of Proposition \ref{constructible}, using Lemma
\ref{Eulermulti}, we can prove that there exists a finite subset
$R$ of $\mo_i\times\mo_j\times\mo_k$ such that
$$
\overline{\mo_i\times\mo_j\times\mo_k}=\bigcup_{(X_i,X_j,X_k)\in
R}\overline{\lr{X_i,X_j,X_k}}
$$
and
$$
F_{\mo_i\mo_j}^{\mo_k[1]}=\sum_{(X_i,X_j,X_k)\in
R}\chi(\overline{\lr{X_i,X_j,X_k}})\cdot
F_{\mo_{X_i}\mo_{X_j}}^{\mo_{X_k}[1]} .
$$
and
$$
F_{\mo_j\mo_k}^{\mo_i[1]}=\sum_{(X_i,X_j,X_k)\in
R}\chi(\overline{\lr{X_j,X_k,X_l}})\cdot
F_{\mo_{X_j}\mo_{X_k}}^{\mo_{X_i}[1]} .
$$
This implies our claim.

If $h_{\mo_1}+h_{\mo_2}+h_{\mo_3}\neq 0$, then any term in $a_M$
vanishes. So we assume $h_{\mo_1}+h_{\mo_2}+h_{\mo_3}=0$, in this
case, $a_M=0$ follows our claim. Now consider $\Delta',$
$$
\hspace{-4cm}\Delta'=[[h_{\mo_1},1_{\mo_2}],1_{\mo}]-[[h_{\mo_1},1_{\mo}],1_{\mo_2}]-[[1_{\mo},1_{\mo_2}],h_{\mo_1}]
$$$$\hspace{-2cm}=(h_{\mo_1}\mid
h_{\mo_2}+h_{\mo})[1_{\mo_2},1_{\mo}]+(h_{\mo_1}\mid
h_{\mo_2}+h_{\mo})[1_{\mo},1_{\mo_2}]=0
$$

Finally according to the definition of the Lie bracket, it is easy
to prove $\Delta''=0$. Thus we  complete the proof of Theorem {\bf
C}.\hfill$\square$
\bigskip
\\
{\bf Proof of Theorem D:}  We claim $( , )\mid_{\fn\times\fn}$ is
non degenerated. For any dimension vector $\bd$ and  $\hat{f}\in
I_{GT}(\bd)$, without loss of generality, we may assume $\mo_i\cap
\mo_j=\emptyset$ for $i\neq j.$ If $\hat{f}\neq 0,$  then there
exists $m_i\neq 0.$ We take any $L\in \mo_i$ and let $\mo_{L[1]}$
be the orbit of $L[1]$, then $(\hat{f}\mid
\hat{1}_{\mo_{L[1]}})=m_i\neq 0$. So it remains to prove that the
bilinear form defined by $(b),(c)$ and $(d)$ is symmetric and
satisfies the condition $(a)$. For $(i)$,
$$
([h_{\bd},\hat{1}_{\mo_1}]\mid \hat{1}_{\mo_2})=(h_{\bd}\mid
h_{\mo_1})(\hat{1}_{\mo_1}\mid \hat{1}_{\mo_2})=(h_{\bd}\mid
h_{\mo_1})\chi(\overline{\mo_1\cap \mo_{2}[1]}) $$ and
$$
(h_{\bd}\mid[\hat{1}_{\mo_1},\hat{1}_{\mo_2}])=(h_{\bd}\mid
h_{\mo_{1}\cap \mo_{2}[1]})\chi(\overline{\mo_1\cap
\mo_{2}[1]})=(h_{\bd}\mid h_{\mo_1})\chi(\overline{\mo_1\cap
\mo_{2}[1]})
$$
For $(ii)$, we only need to prove the following identity:
$$
([\hat{1}_{\mo_{1}},\hat{1}_{\mo_{2}}]_{\mathfrak n}\mid
\hat{1}_{\mo_3})=(\hat{1}_{\mo_{1}}\mid
[\hat{1}_{\mo_{2}},\hat{1}_{\mo_{3}}]_{\mathfrak n})
$$
We have known the Lie bracket  of constructible functions is still
a constructible function. Hence, we can set
$[\hat{1}_{\mo_1},\hat{1}_{\mo_2}]_{\mathfrak n}=\sum
m_{i}\hat{1}_{\mathcal{C}_i}$ and
$[\hat{1}_{\mo_2},\hat{1}_{\mo_3}]_{\mathfrak n}=\sum
n_{j}\hat{1}_{\mathcal{D}_j}$. For any $L\in \mathcal{C}_i$ and
$M\in \mathcal{D}_j$,
$$
m_{i}=F_{\mo_1\mo_2}^{L}-F_{\mo_2\mo_1}^{L} \mbox{ and
}n_{j}=F_{\mo_2\mo_3}^{M}-F_{\mo_3\mo_2}^{M}.
$$
Therefore,
$$
([\hat{1}_{\mo_{1}},\hat{1}_{\mo_{2}}]_{\mathfrak n}\mid
\hat{1}_{\mo_3})=\sum m_{i}\chi(\overline{\mathcal{C}_{i}\cap
\mo_{3}[1]})=F_{\mo_1\mo_2}^{\mo_{3}[1]}-F_{\mo_2\mo_1}^{\mo_{3}[1]}
$$
$$
(\hat{1}_{\mo_{1}}\mid
[\hat{1}_{\mo_{2}},\hat{1}_{\mo_{3}}]_{\mathfrak n})=\sum
n_{j}\chi(\overline{\mathcal{D}_{j}\cap
\mo_{1}[1]})=F_{\mo_2\mo_3}^{\mo_{1}[1]}-F_{\mo_3\mo_2}^{\mo_{1}[1]}
$$
\nd Using our claim (\ref{1649}) in the proof of Theorem {\bf C},
we obtain the identity. \hfill$\square$

\subsection{}

Let $K_0$ be the Grothendieck group of $\md^b(A)$ and
$\mathfrak{h}_{\bbz}=K_0.$  For $\bd\in K_0,$ let
$M_{GT}^{\bbz}(\bd)$ and $I_{GT}^{\bbz}(\bd)$ be the set of all
$\hat{f}\in M(\mp_2(A,\bd))$ and $\hat{f}\in I_{GT}(\bd)$ such that
all values of $\hat{f}$ are in $\bbz$, respectively. We consider the
$\bbz$-spaces $\fn_{\bbz}=\bigoplus_{\bd\in K_0}I_{GT}^{\bbz}(\bd)$
and
$\mathfrak{g}_{\bbz}=\mathfrak{h}_{\bbz}\oplus\mathfrak{n}_{\bbz}.$
Then Theorem {\bf C} and Theorem {\bf D} have the following
$\bbz$-form.

\begin{theorem}\label{xxz1}
\begin{enumerate}
\item[(1)] The $\bbz$-space $\fg_{\bbz}$ with the bracket $[-,-]$
is a Lie algebra. \item[(2)] The symmetric bilinear form $(-|-)$
is invariant over $\fg_{\bbz}$ in the sense
\\ $([x,y]|z)=(x|[y,z])$ for any $x,y\ \mbox{and}\ z\in\fg_{\bbz}.$
\item[(3)] $(-|-)|_{\fn_{\bbz}\times\fn_{\bbz}}$ is
non-degenerated.
\end{enumerate}
\end{theorem}

\bigskip
\section{Realization of generalized Kac-Moody Lie algebra}
\subsection{} As an application, we consider the case  $A=\bbc Q$ where $Q$ is a finite quiver
without oriented cycle.  In this section, we always use $\fg$ to
denote the Lie algebra arising from $\mp_2(A)$ which is defined in
Section 5.2. We will show that the above Lie algebra is a
generalized Kac-Moody Lie algebra and the corresponding symmetric
Kac-Moody algebra is a subalgebra of it, by following the methods
in \cite{SV2001} and \cite{DengXiao2002}. When $Q$ is a tame
quiver, we will demonstrate that the precise structure of the
symmetric affine Lie algebra can be revealed from the derived
category of representations of $Q.$ This phenomenon was first
discovered in \cite{FMV2001} and extended in \cite{LinPeng2002} .

We fix the embedding $\mod A$ in $D^{b}(A)$ by taking any $X\in
\mod A$ as a stalk complex $X^{\blt}=(X_i)$ with $X_0=X$ and
$X_i=0$ for $i\neq 0.$  Hence, we have the induced embedding of
$\mod A$ in $\md_2(A)$ since the functor $F$ in Section 4.1 is
dense (See \cite{PengXiao1997} Corollary 7.1 or
\cite{keller2005}). Moreover, $\md_2(A)\cong \mk_2(\mp(A)).$

Let $\ind\md_2(A)$ be the set of isomorphism classes of all
indecomposable objects in $\md_2(A)$ and $\ind(A)$ be the set of
isomorphism classes of all indecomposable objects in $\mod(A)$. We
know that $\ind\md_2(A)=\ind(A)\cup (\ind(A)[1])$ by
\cite{Happel1988}.

For any $\alpha\in K_0(\mod A),$ let $\bbe(A,\alpha)$ be the
module variety on the dimension vector $\alpha.$   Let
$M_{G}(\alpha)$ be the $\bbc$-space of $G_{\alpha}$-invariant
constructible functions over $\bbe(A,\alpha).$ We can define a
convolution product over  $M_{G}(A)=\bigoplus_{\alpha\in K_0(mod
A)}M_{G}(\alpha)$  such that
$$
1_{\mo_1}*1_{\mo_2}(M)=\chi(\{0\subset M_1\subset M\mid M_1\cong
Y,M/M_1\cong X \mbox{ for some }X\in \mo_1, Y\in\mo_2\}).
$$
for any $M\in \bbe(A,\alpha_1+\alpha_2)$ and  any
$G_{\alpha}$-invariant constructible subsets $\mo_1\subset
\bbe(A,\alpha_1)$ and $\mo_2\subset \bbe(A,\alpha_2)$ (see
\cite{Riedtmann1994} and \cite{lusztig2000}).

A function is called indecomposable if all points in its support
set correspond to  indecomposable $A$-modules. Let $I_{G}(\alpha)$
be the $\bbc$-space of $G_{\alpha}$-invariant indecomposable
constructible functions over $\bbe(A,\alpha).$ We set
$$\fn_{A}=\bigoplus_{\alpha\in K_0(modA)}I_{G}(\alpha)$$
The convolution induces a well-defined Lie bracket over $\fn_{A}$
(see \cite{Riedtmann1994} and \cite{dxx2006}).  For any $f\in
I_{G}(\alpha),$ there is an equivalence class $\hat{f}$ over
$\mp_2(A, \alpha)$ by the embedding $\mathrm{mod} A\hookrightarrow
\mk_2(\mp(A)).$
 It is not difficult to prove that the morphism:
 $$\fn_{A}\rightarrow \fg $$
 $$\hspace{0.2cm}1_{\mo}\mapsto \hat{1}_{\mo}$$ is an injective Lie algebra homomorphism. So $\fn_{A}$
is isomorphic to a Lie subalgebra (denoted by $\fn^+$) of $\fg,$
which deduces a triangular decomposition
$\mathfrak{g}=\fn^{+}\oplus \mathfrak{h}\oplus \fn^{-}$ and
$\fn=\fn^{+}\oplus \fn^{-}$ where $\fn^-:=\fn_{A}[1]$ is generated
by $\hat{1}_{\mo[1]}$ with $1_{\mo}\in\fn_{A}. $ The triangular
decomposition $\fg=\fn^+\oplus \fh \oplus \fn^-$ and the
non-degenerated bilinear form $(-|-)\mid_{\fn_{A}\oplus
\fn_{A}[1]}$ guarantee that $\fg$  is a generalized Kac-Moody Lie
algebra by \cite{Borcherds1988}. The generators and generating
relations can be constructed by the following process.

\subsection{}
First, Let $\mathcal{LE}(A)$ be the Lie subalgebra of
$\mathfrak{g}$ generated by all $\hat{1}_{\mo_{X}}$ which
$\mo_{X}$ is the orbit of $X$ with $X$ exceptional object in
$\md_{2}(A).$ It can be proved that $\mathcal{LE}(A)$ is the
corresponding Kac-Moody Lie algebra (see\cite{PengXiao2000}
Theorem 4.7). Let $\mathcal{LE}(A)^{\pm}=\mathcal{LE}(A)\cap
\fn^{\pm}$ be its Lie subalgebras of positive and negative parts,
respectively. We set $\mathfrak{g}^{+}_{0}(A)=\mathcal{LE}(A)^{+}$
and $\mathfrak{g}^{-}_{0}(A)=\mathcal{LE}(A)^{-}$. We shall
construct the Lie subalgebra $\mathfrak{g}^{\pm}_{m}(A)$ of
$\mathfrak{n}^{\pm}(A)$ and the Lie subalgebra
$\mathfrak{g}_{m}(A)$ of $\mathfrak{g}(A)$ for $m\geq 0$.

Let $m\geq 1$ and suppose $\mathfrak{g}^{\pm}_{m-1}(A)$ and
$\mathfrak{g}_{m-1}(A)=\fg^-_{m-1}\oplus\fh\oplus\fg^+_{m-1}$ have
been constructed. We let $\pi_{m}\in \mathbb{N}[I]$ have the
smallest trace such that $\mathfrak{g}^{+}_{m-1}(A)_{\pi_m}\neq
\mathfrak{n}^{+}(A)_{\pi_m}$. Then we define
$$
L^{+}_{\pi_m}=\{x^{+}\in \mathfrak{n}^{+}(A)_{\pi_m}\mid
(x^{+}|\mathfrak{g}^{-}_{m-1}(A)_{\pi_m})=0\}
$$
and
$$
L^{-}_{\pi_m}=\{y^{-}\in \mathfrak{n}^{-}(A)_{\pi_m}\mid
(\mathfrak{g}^{+}_{m-1}(A)_{\pi_m}|y^{-})=0\}
$$
We now denote by $\mathfrak{g}^{\pm}_{m}(A)$ the Lie subalgebra of
$\mathfrak{n}^{\pm}(A)$ generated by $\mathfrak{g}^{\pm}_{m-1}(A)$
and $L^{\pm}_{\pi_m}$, respectively. Set
$\mathfrak{g}(A)_{m}=\mathfrak{g}^{+}_{m}(A)\oplus
\mathfrak{h}\oplus \mathfrak{g}^{-}_{m}(A)$. As a conclusion, we
obtain the chains of Lie subalgebra of $\mathfrak{n}^{\pm}(A)$ and
of $\mathfrak{g}(A)$:
$$
\mathfrak{g}^{\pm}(A)_{0}\subset \mathfrak{g}^{\pm}(A)_{1}\subset
\cdots \subset \mathfrak{g}^{\pm}(A)_{m}\subset \cdots \subset
\mathfrak{n}^{\pm},
$$
$$
\mathfrak{g}(A)_{0}\subset \mathfrak{g}(A)_{1}\subset \cdots \subset
\mathfrak{g}(A)_{m}\subset \cdots \subset \mathfrak{g}.
$$
For $m\geq 1$, let
$\eta_{m}=\mathrm{dim}_{\mathbb{C}}L^{+}_{\pi_m}=\mathrm{dim}_{\mathbb{C}}L^{-}_{\pi_m}$.
All $\pi_m$ lie in the fundamental set, that is, $(\pi_m,i)\leq 0$
for $i\in I$ and  $m\geq 1$ (see \cite{DengXiao2002}). Depending
on non-degeneracy of $(-|-)\mid_{\fn^{+}\times\fn^{-}},$ we can
construct a basis $\{E_{p}(m)\mid 1\leq p \leq \eta_{m}\}$ of
$L^{+}_{\pi_m}$  and a basis $\{F_{p}(m)\mid 1\leq p\leq
\eta_{m}\}$ of $L^-_{\pi_m}$ such that
$$
(E_{p}(m),F_{q}(n))=\delta_{pq}\delta_{mn}.
$$
Now we set $e_{i}=E_{i}(0)=\hat{1}_{\mo_{S_i}}$ and
$f_{i}=F_{i}(0)=\hat{1}_{\mo_{S_{i}[1]}}$ for $i\in I$, where
$S_i$ is the simple $A$-module at the vertex $i$. Let
$e_{j}=E_{p}(m)$ and $f_{j}=F_{p}(m)$ for $j=(\pi_{m},p)\in J$
where $J=\{(\pi_{m},p)\mid 1\leq p \leq \eta_{m} \}$. If
$j=(\pi_m,p)\in J$ we denote $h_j=h_{\pi_m}.$ Following the above
construction, $\mathfrak{g}$ is the Lie algebra generated by
$e_{i}, f_{i}$ and $\mathfrak{h}$ for $i\in I\cup J$. We define
$\omega(e_i)=f_i$ for $i\in I\cup J$ and $\omega(h_i)=-h_i$ for
$i\in I.$ The property of non-degeneracy of
$(-|-)\mid_{\fn^{+}\times\fn^{-}}$ implies that
$\omega:\fg\rightarrow\fg$ is an involution.

The triangular decomposition $\fg=\fn^-\oplus\fh\oplus\fn^+,$ the
non-degeneracy and invariance of bilinear form $(-|-)$ guarantee
the generators satisfy the generating relations of the generalized
Kac-Moody Lie algebra (see \cite{SV2001} and \cite{DengXiao2002}).
Hence, we have the following theorem by \cite{Borcherds1988}:
\begin{theorem} Let $A=\bbc Q.$ The Lie algebra
$\mathfrak{g}$ is a generalized Kac-Moody Lie algebra with the
following generating relations:
\begin{itemize}
\item[(1)] $[h_i,h_j]=0$ for any $i,j\in I$; \item[(2)] $[h_i,
e_j]=(i,j)e_j$ for any $i\in I$ and $j\in I\cup J$;
\item[(3)]$[h_i, f_j]=-(i,j)f_j$ for any $i\in I$ and $j\in I\cup
J$;\item[(4)] $[e_i,f_j]=\delta_{ij}h_j$ for any $i,j\in I\cup J$;
\item[(5)]$(ade_{i})^{1-a_{ij}}e_{j}=(adf_{i})^{1-a_{ij}}f_{j}=0$
for any $i\in I$, $j\in I\cup J$ and $i\neq j$;
\item[(6)]$[e_i,e_j]=[f_i,f_j]=0$ if $a_{ij}:=(i,j)=0$ for $i,j\in
I\cup J$.
\end{itemize}
Moreover, the Lie subalgebra $\fg_0(A)$ generated by $e_i,f_i$ and
$h_i$ for $i\in I$ is the derived Kac-Moody Lie algebra with
Cartan datum $(I,(-,-))$, where $(-,-)$ is the symmetric Euler
bilinear form of $Q.$
\end{theorem}

\subsection{}
In the following part, we will suppose $Q$ is a Dynkin or tame
quiver. An explicit construction of the positive part of affine
Kac-Moody algebra has been given in  \cite{FMV2001} via Hall
algebra approach. We will verify that, if we apply our result to
the case the derived category of representations of a tame quiver,
the main result in \cite{FMV2001} can be extended easily to give a
global realization of the affine Kac-Moody algebra. We remark here
that the similar result has been obtained by \cite{LinPeng2002} in
a different way.

\nd{6.3.1} As showed in \cite{GM1996} and \cite{XZZ2005},
BGP-functor( see \cite{BGP1973}) can be defined over root category
$\md_2(A),$ denoted by:
$$
\md_2(A)\xrightarrow{H_{m}(\mathcal{S}^{+}_a)(H_{m}(\mathcal{S}^{-}_a))}\md_2(\sigma_aA)
$$
where $a$ is a source (sink) of quiver $Q$ and
$\sigma_aA=\bbc(\sigma_aQ).$ In particular,
$$H_{m}(\mathcal{S}^{+}_a)(S_a)=S'_a[1]$$ where $S_a$ and $S'_a$ are
the simple modules at the vertex $a$ in $\mod A$ and $\mod
\sigma_a A$, respectively.  The explicit construction as follows
(see Ex 4.6 in \cite{GM1996}). Let $M^{\bullet}=(M_n,\pt_n)\in
\md_2(A)$ where $M_n=(\bbc^{\alpha_n},x_n)$ with dimension vector
$\alpha_i$ for $n\in \bbz$. Recall
$\bbc^{\alpha_n}=\bigoplus_{i\in Q_0}\bbc^{d_i(\alpha_n)}$ and
$x_n=(x_{n,\alpha})_{\alpha\in Q_1}$ with $x_{n,\alpha}\in
\mathrm{Hom}(\bbc^{d_{s(\alpha)}(\alpha_n)},
\bbc^{d_{t(\alpha)}(\alpha_n)})$. Of course, $M_n=M_{n+2}$. Define
$H_m(S^+_i)(M^{\bullet})=(N_n, d_n)\in \md_2(\sigma_i A)$ as
follows. Assume that $N_n=(\bbc^{\beta_n}, y_n)$ with dimension
vector $\beta_n$. Then
$$
\bbc^{d_i(\beta_n)}=\bbc^{d_i(\alpha_n)} \mbox{ for any } i\neq a,
$$
and
$$
\bbc^{d_a(\beta_n)}=(\bigoplus_{(ab)\in
Q_1}\bbc^{d_b(\alpha_n)})\oplus \bbc^{d_a(\alpha_{n+1})}
$$
and for $\beta\in Q_1,$
$$y_{n,\beta}=\left\{\begin{array}{cc}
x_{n,\beta},& {\rm if}\,\  t(\beta)\neq a,\\
{\rm the\ embedding},&{\rm otherwise.}
\end{array}
\right.$$ Recall that $\mathcal{LE}(Q)$ is Lie subalgebra of
$\mathfrak{g}$ generated by the character functions of exceptional
objects in the root category. This functor induces the
homeomorphism between $\mp_2(A)$ and $\mp_2(\sigma_iA)$ so that we
have the following isomorphism between the corresponding Lie
algebras:
$$
\tilde{\varphi}_{i}:\mathcal{LE}(A)\rightarrow
\mathcal{LE}(\sigma_{i}A)
$$
By Theorem 6.1 we have the following canonical isomorphism:
$$
\phi: \mathfrak{g}(C)\rightarrow \mathcal{LE}(A)
$$  fixing the generators $e_i$ and $f_i$ for $i\in I,$
where $C=(I,(-,-))$ is the Cartan datum of  $Q$, $\mathfrak{g}(C)$
is the derived Kac-Moody Lie algebra corresponding to $C$.
Moreover, we have the following commutative diagram (see
\cite{XZZ2005}):
$$
\xymatrix{\mathfrak{g}(C) \ar[r]^{\tilde{s}_i} \ar[d]^{\phi}&
\mathfrak{g}(C)
\ar[d]^{\phi} \\
\mathcal{LE}(A) \ar[r]^{\tilde{\varphi}_{i}} &
\mathcal{LE}(\sigma_{i}A)}
$$
where
$\tilde{s}_i:=\mathrm{exp}(ade_i)\mathrm{exp}(ad(-f_i))\mathrm{exp}(ade_i)$
is the isomorphism map defined in \cite{kac1990}.

\bigskip
\nd{6.3.2} \hspace{0.3cm} In order to give the concrete expression
of Chevalley basis in $\fg(C)$ by the above isomorphism $\phi$ in
$\mathcal{LE}(A),$ we will consider Lie algebra based on the Euler
cocycle. Let $Q$ is a tame quiver and $R$ be the root system of
the quiver $Q$. Let $\delta$ be the minimal imaginary root of
$R_{+}.$ We denote by $\mathfrak{g}^{\epsilon} (C)$ the following
$R$-graded $\mathbb{C}$-linear space:
\begin{equation}\nonumber
\mathfrak{g}^{\epsilon} (C) = \bigoplus_{\alpha \in R\cup{\{0\}}}
\mathfrak{n}_{\alpha}^{\epsilon} (C) ,
\end{equation}
such that
\begin{itemize}
\item for each  real root $\alpha$ a one-dimensional
$\mathbb{C}$-linear space $\mathfrak{n}_{\alpha}^{\epsilon} (C) =
\mathbb{C} \Tilde{e}_{\alpha} $ with generator
$\Tilde{e}_{\alpha}$, \item for each  imaginary root $n \delta$ a
$\mathbb{C}$-linear space $\mathfrak{n}_{n \delta}^{\epsilon} (C)
= \mathbb{C} [I] / \mathbb{C} \delta$, where we consider $\delta$
as an element of $R_+ \subset \mathbb{Z} [I] \subset \mathbb{C}
[I]$. For $h \in \mathbb{C} [I]$ we denote by $h (n)$ the image of
$h$ under the natural projection map $\mathbb{C} [I] \rightarrow
\mathfrak{n}_{n \delta}^{\epsilon} (C)$.
\end{itemize}

The space $\fg^{\varepsilon}(Q)$ is equipped with the bilinear
bracket

\begin{equation}
\begin{split}
[\Tilde{e}_\alpha, \Tilde{e}_\beta] &=
\begin{cases}
\epsilon (\alpha , \beta ) \Tilde{e}_{\alpha + \beta} &
\text{ if } \alpha + \beta \in R^{\re} , \\
\epsilon (\alpha , \beta ) \alpha (k) &
\text{ if } \alpha + \beta = k \delta , \\
h_{\alpha} & \text{ if } \alpha =-\beta ,\\
0 & \text{ otherwise},
\end{cases}
\\
[h (n), \Tilde{e}_{\alpha} ] = -[\Tilde{e}_{\alpha} , h(n)] &=
\epsilon (n \delta , \alpha ) (h, \alpha) \Tilde{e}_{\alpha + n
\delta} ,
\\
[h(n),h'(-n)]] &=n(h,h')c,\  c \text{ is the center element}
\\
[h (n), h (m) ] &= 0 ,\text{ if }n+m\neq 0
\end{split}\end{equation}
where $h \in \mathbb{C} [I]$, $\epsilon$ is the Euler cocycle .
The following is well-known.

\begin{theorem}
The $\mathbb{C}$-linear space $\mathfrak{g}^{\epsilon} (Q)$
equipped with the above bracket  is isomorphic to the derived
affine Kac-Moody algebra of type $C$.
\end{theorem}

Next, we describe the structure of $\md_2(A)$ (See
\cite{Happel1988}) for $A=\bbc Q$ of tame quiver $Q$. We know
$\ind\md_2(A)=\ind A\cup\ind A[1].$ The AR-quiver of $\md_2(A)$
consists of two components of form $\bbz Q$ and two families of
orthogonal  tubes indexed by $\mathbb{P}^{1}(\bbc).$  We denote
the two families of tubes by $\mathbb{T}_{Q}=\{T_{z}\mid z\in
\mathbb{P}^{1}(\bbc)\}$ and by $N_{z}$ the period of $T_z$ and
$\mathbb{T}_{Q}[1]=\{T_{z}[1]\mid z\in \mathbb{P}^{1}(\bbc)\}.$
Let $L\subset \mathbb{P}^{1}(\bbc)$ be such that the tube $T_z$ is
of period $N_z>1$ for any $z\in L.$ Then $|L|\leq 3.$  By
\cite{DlabRingel}, we may fix
an embedding: $\pi:\md_2(K)\rightarrow \md_2(A)$ such that
$\pi(\mathbb{T}_{K})\subset\mathbb{T}_{Q}$ and
$\pi(\mathbb{T}_{K}[1])\subset\mathbb{T}_{Q}[1]$ for the Kronecker
quiver $K.$ It induces the Hall map $\pi_{*}$ from $M_{GT}(K)$ to
$M_{GT}(A)$.

Let $\hat{1}_{(n,n)}\in I_{GT}(K, (n,n))$ be the equivalence class
over $\mp_2(K, (n,n))$ of the characteristic function  on the
constructible subset of all indecomposable objects in
$\mathrm{mod} K$ with dimension vector $(n,n)$ for $n\in\bbz.$ We
denote by $E_0(n)$ the image of $\hat{1}_{(n,n)}$ under $\pi_{*}
.$

Let $M_{i,1,z}$ be the regular simple objects in $T_{z}$ and
$M_{i,l,z}$ be the indecomposable objects in $T_{z}$ with regular
length $l$ for $i\in \bbz/N_z \bbz$ such that the regular socle of
$M_{i,l,z}$ is $M_{i,1,z}.$  Set $M_{i,-l,z}=M_{i,l,z}[1]$ for
$l>0.$ Let $E_{i,l,z}:=1_{\mo_{M_{i,l,z}}}$ be the characteristic
function of the orbit of $M_{i,l,z}.$  If $|l|\not\equiv
0(modN_z),$ then \underline{\textbf{dim}}$M_{i,l,z}\in R^{re}$ and
$E_{i,l,z}\in\mathcal{LE}(A).$ If $l=nN_z,$ then
\underline{\textbf{dim}}$M_{i,l,z}=n\delta$ and
$E_{i+1,l,z}-E_{i,l,z}\in \mathcal{LE}(A)$ (See \cite{FMV2001} and
\cite{PX1996}).

We define the function $\xi:R\rightarrow {\pm1}$ by
$$
\xi(\alpha)=(-1)^{(1+\mathrm{dim_{\bbc}End}M)}$$ for
indecomposable object $M$
 with $\underline{\textbf{dim}}M=\alpha$.

For any function $f_{\alpha}\in I_{GT}(\md_2(A)),$ set
$\tilde{f}_{\alpha}=\xi(\alpha)f_{\alpha}.$

Now we are ready to give the generalization of the main theorem in
\cite{FMV2001}. By using the reflection functor defined for
$\md_2(A)$ and following the method in \cite{FMV2001} we
eventually obtain the following result.
\begin{theorem}\label{AffineTheorem}
Let $Q$ be any tame quiver without oriented cycles. The following
map completely describes an isomorphism
$\fg^{\epsilon}(C)\rightarrow \mathcal{LE}(A)$:
\begin{equation}\nonumber
\begin{split}
\Xi^{\epsilon} (\Tilde{e}_{\alpha}) = \Tilde{E}_{\alpha}
\text{ for any $\alpha\in R^{\re}$ } , \\
\Xi^{\epsilon} (\alpha_{i,z} (n)) =
\Tilde{E}_{i,nN_z,z} - \Tilde{E}_{i+1,nN_z,z} ,\\
\hspace{-5cm}\Xi^{\epsilon} (\alpha_{0} (n)) = - \Tilde{E}_{0} (n) \\
\Xi^{\epsilon} (\beta) = h_{\beta}, \Xi^{\epsilon} (c) =
h_{\delta}.
\end{split}
\end{equation}
for any $n\in\bbz$ and $\beta\in \bbz[I]$.
\end{theorem}

Obviously, the set:
$$
\{\Tilde{E}_{\alpha}\}_{\alpha\in R^{\re}}\cup
\{\Tilde{E}_{i,nN_z,z} -
\Tilde{E}_{i+1,nN_z,z}\}_{i\in\bbz/N_{z}\bbz,z\in L,n\in\bbz}\cup
\{\Tilde{E}_{0} (n)\}_{n\in\bbz}\cup\{h_i\}_{i\in I}
$$
provides a $\bbz$-basis of $\mathcal{LE}(A).$

\end{document}